\documentclass[final,3p,10pt]{elsarticle}

\usepackage{graphicx}
\usepackage{amssymb}
\usepackage{amsthm}

\usepackage{amsmath}
\usepackage{amsfonts}
\usepackage{caption}
\usepackage{subcaption}
\usepackage{algorithmicx}
\usepackage[ruled]{algorithm}
\usepackage{algpseudocode}
\usepackage{varioref}

\theoremstyle{definition}
\newtheorem{defin}{Definition}
\newtheorem{propi}{Theorem}
\newtheorem{hyp}{Assumption}

\newtheorem{rem}{Remark}

\labelformat{defin}{def.~#1}
\labelformat{propi}{th.~#1}
\labelformat{hyp}{ass.~#1}
\labelformat{lem}{lem.~#1}
\labelformat{rem}{rem.~#1}
\labelformat{ex}{ex.~#1}

\newcommand{\dif}{{\mathrm d}}

\newcommand{\egdef}{{\, \stackrel{\text{def}}{=} \,}}
\newcommand{\egdist}{{\, \stackrel{\text{dist}}{=} \,}}

\newcommand{\ensM}{{\mathbb{M}}} 
\newcommand{\ensS}{{\mathbb{S}}} 
\newcommand{\ensU}{{\mathbb{U}}} 

\newcommand{\tribM}{{\mathsf{M}}} 
\newcommand{\tribS}{{\mathsf{S}}}

\newcommand{\mesM}{{\mathcal{M}}} 
\newcommand{\mesS}{{\mathcal{S}}} 
 
\newcommand{\mesT}{{\mathcal{T}}} 

\newcommand{\N}{{\mathbb{N}}}
\newcommand{\R}{{\mathbb{R}}}

\newcommand{\Rt}{{\mathbb{R}^3}}

\newcommand{\thx}{{\theta_{xy}}}
\newcommand{\thz}{{\theta_{z}}}

\begin{document}

\begin{frontmatter}


\title{Stochastic simulation of urban environments \\ Application to Path-loss in wireless systems}

\author[telecom]{Thomas Courtat}
\author[telecom]{Laurent Decreusefond}
\author[telecom]{Phillipe Martins}
\address[telecom]{T\'{e}l\'{e}com ParisTech, 23 avenue d'Italie, 75013 Paris, France}

\begin{abstract}
We are interested in the assessment of electromagnetic Path-Loss in complex environments. The Path-loss is the attenuation function $P$ of the electromagnetic power at a distance $d$ of an antenna. In free-space, $P(d) \propto 1/d^2$, in complex environments like cities, wave trajectory is altered by successive reflections and absorptions, the path-loss is not theoretically known and engineering rules postulate that $P(d) \simeq 1/d^{\gamma}, \, \gamma>2$.\\
We place in a stochastic geometry context to answer the problem statistically. We present random models of 3D-city. These models reproduce main real cities' features, can be calibrated with simple mean formulae and can be fast simulated. 
For collections of random cities with the same mean morphology, we estimate by Monte-Carlo ray tracing techniques their attenuation maps. By averaging these maps, we show that the power expectancy actually follows a function $\sim 1/d^{\gamma}$ with $\gamma$ depending on the environment morphology.
\end{abstract}


\end{frontmatter}


\section{Introduction \label{sec:introduction}}
In this paper, we seek out to analyze quantitatively  the power field generate by an electromagnetic source over a complex urban environment.
 \\
In the ideal case of free space propagation, the power  $\mathbf{P}(d)$ at a point only depens on its distance $d$ to the source: $\mathbf{P}(d) \propto 1/ d^2$.
In a urban environment, the power field in a point is constituted by a multitude of paths created by successive reflections of source waves on the buildings of the city. Qualitative and empirical engineering rules are applied to assess the power's order of magnitude\cite{Sarkar2010}. The formula that is commonly admitted is $\mathbf{P}(d) = \dfrac{c}{d^{\alpha}}\cdot H $ with $2 < \alpha \leq 4$ depending on the specificity of the environment under consideration and $H$ a random variable following a Rayleigh or log-normale law, independent to the distance, to account the interference fading and the shadowing (zones that are partially or totally darkened by elements of the landscape).
\\
The synthetic knowledge of this exponent is of main interest in planning and sizing of wireless telecommunication networks.
\\
We present here a statistical approach to solve this problem. Rather than trying to determine the precise attenuation map for a particular city, we build representative stochastic models for cities and explore the statistical properties of the random attenuation map in each model.
\\
 From a deterministic point of view, the regular square lattice, the "Manhattan" model is often used \cite{Sun2005} to represent the street network. But this model very particular geometrical features that do not give a full account of the variability of morphologies that are observed between real cities.
\\
Stochastic geometry \cite{Stoyan1995,Lieshout2000}, allows proposing random models of planar network \cite{Baccelli1995,Baccelli2009a,Kendall2010} and to calculate statistical indicators to answer various planning problems. We can for instance cite the coverage probability of a cellular network with or without user displacement \cite{Morlot2012} or the length distribution in a cooper or optical network \cite{Gloaguen2009}. 
\\
The outline of stochastic approach is to describe in one hand the city (yet restricted to its street network in the literature) by a small number of parameters (as the street density, the mean number of intersections by $\text{km}^2$...) and in the other hand the structure of the telecommunication network by another small set of parameters (the spatial density of routers, of users...). Theoretical models of networks (Poisson Line, Poisson Vorono\"{i}, or Poisson Delaunay \cite{Gloaguen2002}) and spatial point distribution models (Poisson Processes, Cox Processes \cite{Baccelli1995,Gloaguen2005,Morlot2012}) allow simulating collections of cases with the same statistical properties ; the statistical answer to the problem under consideration is then obtained by evaluating the expectancy of relevant functional over these collections. In some problems or limit cases, this evaluation can be made from analytical calculation \cite{Voss2010} or be rewritten as real valued function integrals one just has to compute numerically \cite{Morlot2012}. 
\\
As for Path-loss function problem, a ray tracing approach \cite{Catedra1998} has been presented in \cite{Marano2005}. The retained city model is a square lattice with percolation: each site is occupied by a building with a given probability. The electromagnetic wave propagates from a punctual source in the horizontal plane by reflecting on buildings. One of their key point to estimate power field is that as it reflects a ray performs a recurrent random walk. Yet if a random walk is recurrent in the plane, this phenomenon does not happen in dimension three. We thus can expect different results if we take the three dimensional structure of the city into account.

The first step of this article presents three dimensional random models for cities, both tractable and realistic. They mimic main features of cities: facades alignment along streets, organization of buildings into blocks, variability in the topology of street intersections and of the network anisotropy \cite{Courtat2011}. We proceed into two parts: at first we generate street axis from random tessellations (Sec.\ref{sec:tessellation}), secondly we add random buildings in the cells of the tessellation (Sec.\ref{sec:batiments}).
\\
The second step is the probabilistic modelling of propagation and the realization of an efficient Monte-Carlo ray tracing algorithm to evaluate the power field generated by an antenna placed in height (Sec.\ref{sec:propagation}).
\\
From this point, we can simulated classes of statistically equivalent cities and compute the power field also called attenuation map and estimate the expectancy of the received power at a distance $d$ of an antenna (Sec.\ref{sec:resultats}).
\\
The algorithms are efficiently implemented under our \textit{GeoStat} framework. Every step of the simulations is done with vector objects: to avoid aliasing effects, buildings are collections of polytopes and the maps are never rasterized. Algorithms we present are optimized, the main point of algorithms are present through the article body and their technical refinement are sent back to annexes (\ref{ann:algo1}). Their need in memory and running time are discussed in (\ref{ann:temps}).


\section{Street Models \label{sec:tessellation}}
\subsection{General definitions}
\begin{defin}[Directed line]
In the Euclidian plane, a directed line $L$ is defined by its origin $o_L$ and its direction vector $u_L$.
\\
The set of lines in the plane is written $\mathbb{H}$
\end{defin}
\begin{defin}[Line]
If $L$ is a directed line then $\vert L \vert = \lbrace o_L + \lambda \vec u_L, \quad \lambda \in \mathbb{R} \rbrace $ is a subset of $\mathbb{R}^2$ called a line.
\\ 
The set of lines is written $\vert \mathbb{H} \vert $.
\end{defin}
A line $L$ is unequivocally defined by $\pi_L(o)$ the projection of the origin $0$ of the plane over $L$ i.e by a signed distance to the origin and a geometric angle with the $x-$axis. 
Thus $\vert \mathbb{H} \vert$ is seen as the cylinder $\mathbb{R} \times [0,\pi[$.
\\
A directed line $L$ divides  $\mathbb{R}^2$ into two disconnected parts: $L_+ = \{ x \in \mathbb{R}^2, \det(u_l,  (x - o_l) ) >  0 \}$ and $L_- = \{ x \in \mathbb{R}^2, \det(u_l,  (x - o_l) ) <  0 \}$. 
\\
From this remark, we adopt the following general recursive definition of polygons: 
\begin{defin}[Polygons]
The set of polygons $\mathbb{P}$ is the subset of convex bodies $\mathbb{K}$ defined by $ \mathbb{P} = \bigcup_n \mathbb{P}_n$ with $\mathbb{P}_0 = \lbrace \mathbb{R}^2 \rbrace$ and 
$\forall n>0, \quad \mathbb{P}_n = \lbrace P \in \mathbb{K}, \quad \exists Q \in \mathbb{P}_{n-1} \text{ and } L \in \mathbb{H}, \quad P = Q \cap L_+ \text{ or } P = Q \cap L_-\rbrace$

\end{defin}

\begin{defin}[Tessellation and cells]
A (straight convex) tessellation $\mathbf{\Xi}$ is a countable family of convex polygons $\{C_i \}_{i \in \mathbb{N}}$ (the cells) partitioning $\mathbb{R}^2$ and whose interiors do not intersect. 
\end{defin}
For instance, if $X =\lbrace X_i \rbrace$ is a countable family of points in the plane, to each point  $X_i$ is associated is Vorono\"{i} zone $V(X_i || X)$. Then $\{V(X_i || X) \}_i$ is a tessellation. 

To the border $\partial \mathbf{\Xi} = \cup \partial C_i$ of a tessellation is associated a straight planar graph denoted $\mathcal{G}(\partial \mathbf{\Xi})$ whose vertices are the cells' vertices and edges are cell's. 
\begin{defin}[Axis]
Two edges of  $\mathcal{G}(\partial \mathbf{\Xi})$ are aligned if they are equal or intersect and make a 0 angle.
\\   
The alignment relationship is symmetric.
\\ 
The transitive closure of the alignement relationship generates a partition of $\partial \mathbf{\Xi}$ deprived of the $\mathcal{G}(\partial \mathbf{\Xi}))$'s vertices into sets composed of one segment, one line or one half-lines. These sets are called axis. 
\end{defin}

\subsection{Elements of random geometry}
\begin{defin}[General Poisson Point Process]
If $\Omega$ is a Polish Space equipped with a Borelian measure $\Lambda$, there exists a countable point process $X = \sum X_i \in\Omega$ such that: 
\begin{enumerate}
\item If $K_1$ and $K_2$ are two measurable sets of $\Omega$ that do not intersect then $X \cap K_1$ and $X \cap K_2$ are independent point processes. 
\item If $K$ is a compact of $\Omega$, the number of points of $X$ in $K$ is a random variable with a Poisson distribution: $\sharp X \cap K \sim \mathcal{P}(\Lambda(K))$
\end{enumerate}
\end{defin} 
For instance in $\mathbb{R}^n$ with $\Lambda(\cdot) = \lambda \cdot \mu(\cdot)$, $\mu$ being the Lebesgue's measure and $\lambda$ a positive real number we obtain the "stationary" Poisson Point Process of intensity $\lambda$

\begin{defin}[Poisson Line (PL)]
The Poisson Line Process (PL) $\partial \Xi$ of intensity measure $\lambda$ and anisotropy probability measure $\mathcal{R}$ is the Poisson Point Process on the cylinder $\mathbb{R} \times [0,\pi[ = \vert \mathbb{H} \vert $ of intensity $\mathcal{L}(.)=\lambda\mu(.) \otimes \mathcal{R}(.)$, $\lambda >0$. 
\end{defin}
\begin{propi}[Random Line]
If $W$ is a compact set in the plane then $\mathcal{L}_W(.)= \mathcal{L}(\{ h \in \mathcal{H}, h \cap W \neq \emptyset \}) < \infty$. \\
A Uniform Random Line in $W$ is a random variable on $|\mathbb{H}|$ with distribution $\mathbb{U}_W= \dfrac{\mathcal{L}_W(.)}{\mathcal{L}_W(\vert \mathbb{H} \vert)}$
\end{propi}
\begin{propi}[Random Tessellation]
Equipped with the $\sigma$-algebra generated by sets $\{\mathbf{\Xi}, \partial\mathbf{\Xi} \cap K = \emptyset   \}$ where $K$ is any compact subset of $\mathbb{R}^2$, the set of tessellation is measurable \cite{Kendall2010}. 
\\
Let thus $\mathbf{\Xi}$ be a random tessellation. Its intensity function is the measure defined by $L_A(W) \egdef \mathbb{E}(\mu_1(\Xi \cap W))$ for any compact $W$.   
\end{propi}
A random tessellation may be deducted from a point a random point process. For instance, $X = \sum X_i$ if is a Poisson Point Process in the plane, then $\{V(X_i || X) \}_i$ is a random tessellation named Poisson Vorono\"{i} Tessellation \cite{Stoyan1995}. 

\begin{defin}
A random tessellation $\Xi$ is stationary if $\Xi + \vec{v} \egdist \Xi, \quad \forall \vec{v} \in \mathbb{R}^2$.  
\end{defin}
\begin{defin}
A random tessellation is locally finite if $L_A(W) < \infty, \quad \forall \: W \text{ compact}$.
\end{defin}

\subsection{Poisson Line Tessellation and Crack STIT Tessellation}
\begin{figure}[h!]
\begin{center}
\begin{subfigure}[b]{0.3\textwidth}
                \centering
                \includegraphics[width =\textwidth]{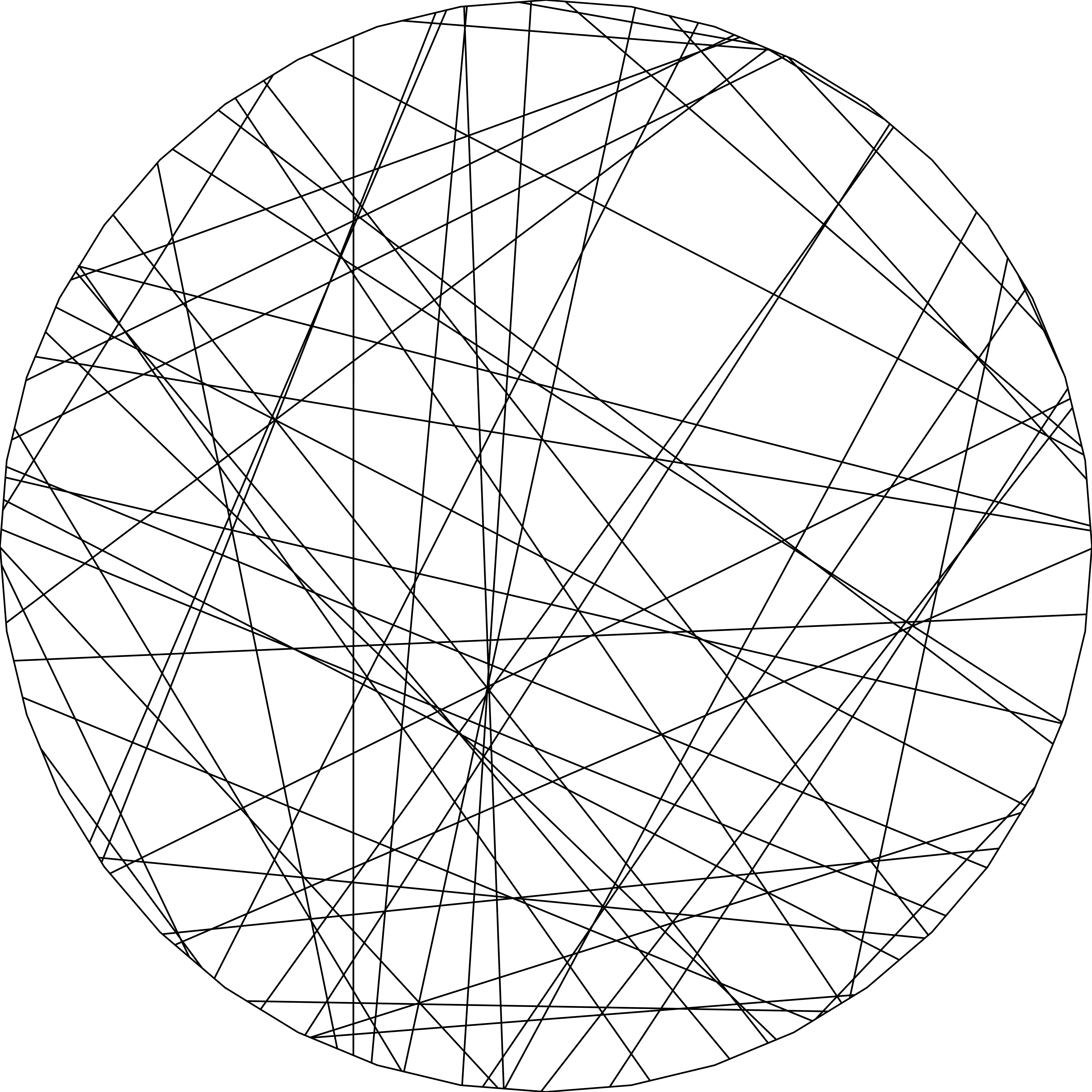}
                
      \end{subfigure}
\begin{subfigure}[b]{0.3\textwidth}
                \centering
                \includegraphics[width =\textwidth]{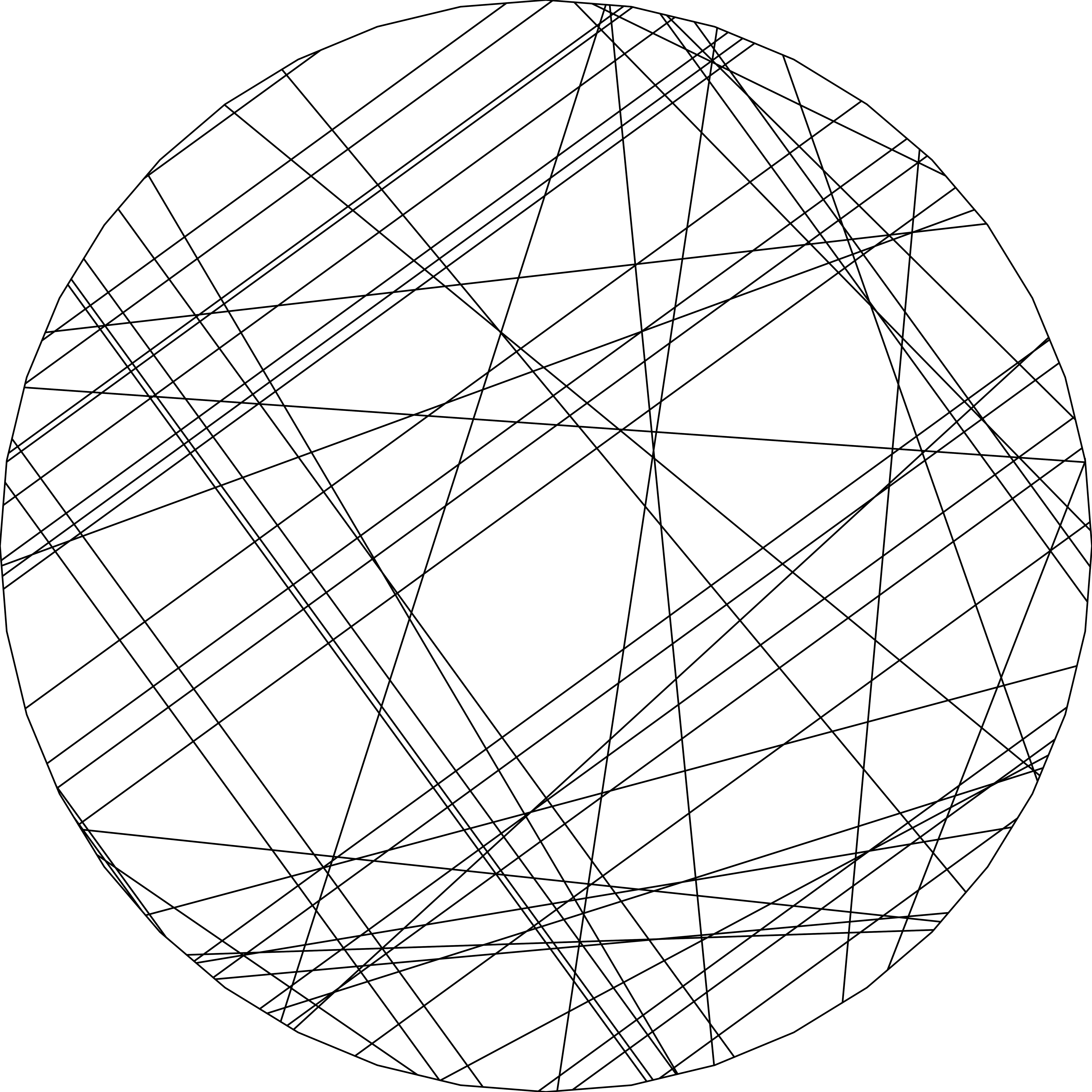}
                
      \end{subfigure}
      \begin{subfigure}[b]{0.3\textwidth}
                \centering
                \includegraphics[width =\textwidth]{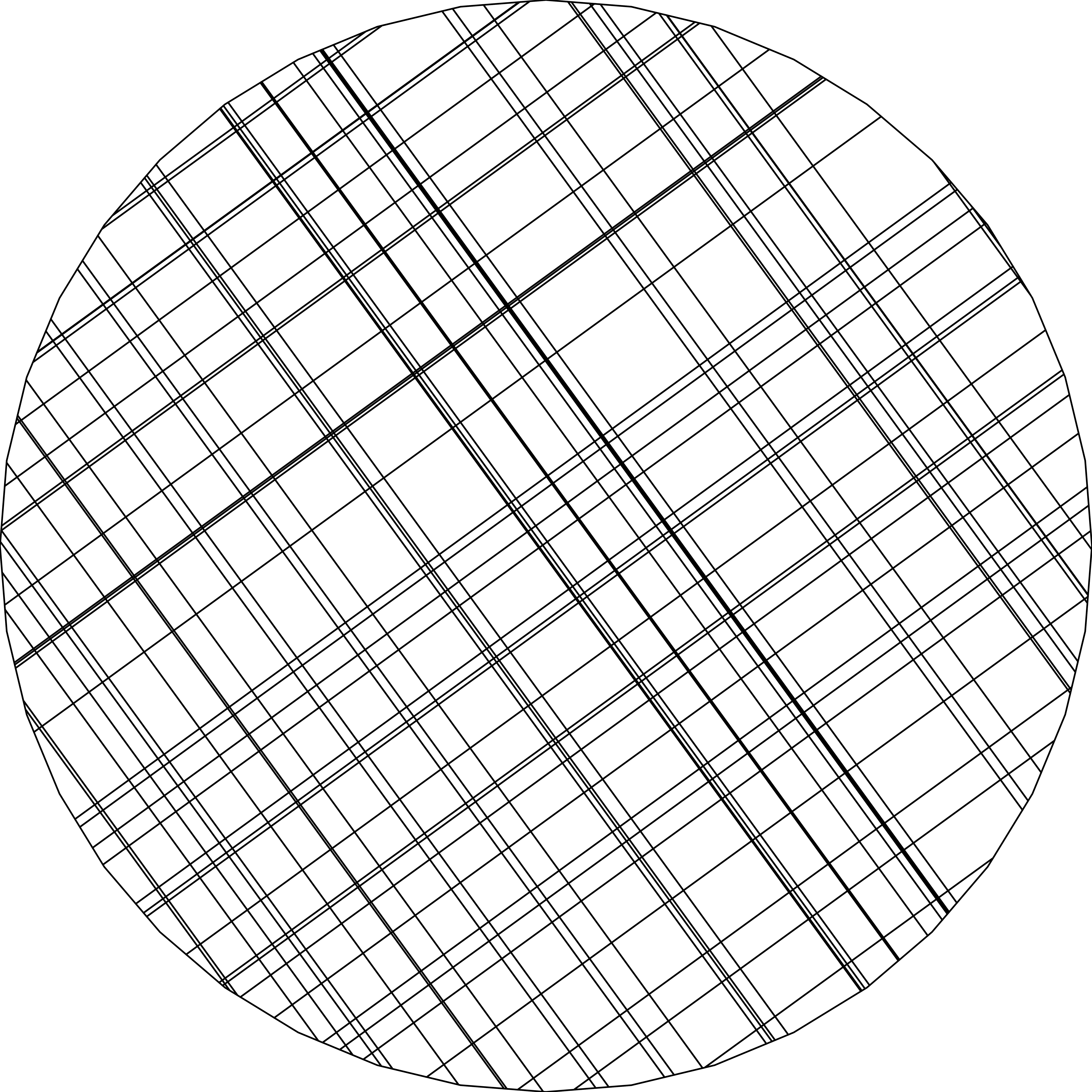}
                
      \end{subfigure}
      \begin{subfigure}[b]{0.3\textwidth}
                \centering
                \includegraphics[width =\textwidth]{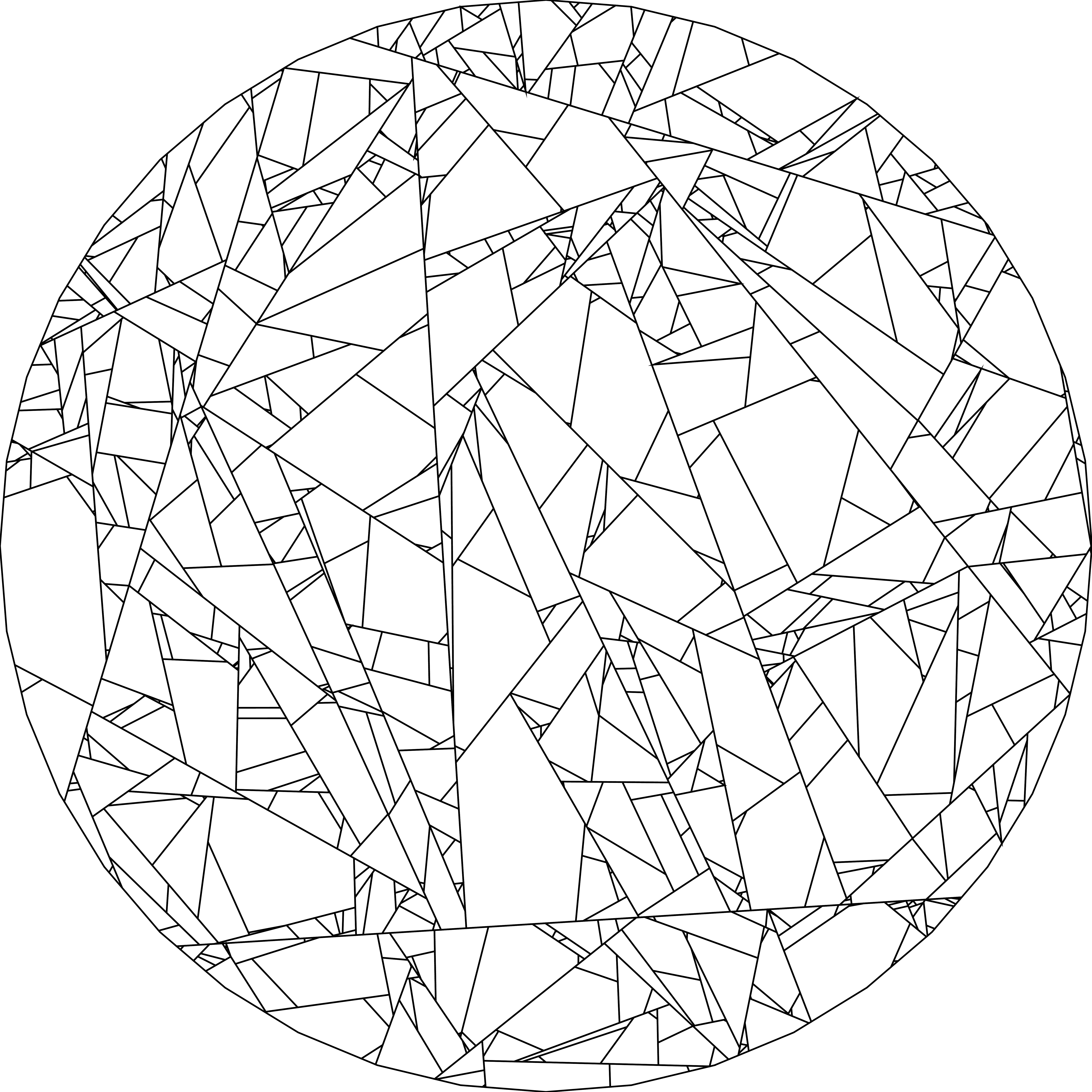}
               
      \end{subfigure}
      \begin{subfigure}[b]{0.3\textwidth}
                \centering
                \includegraphics[width =\textwidth]{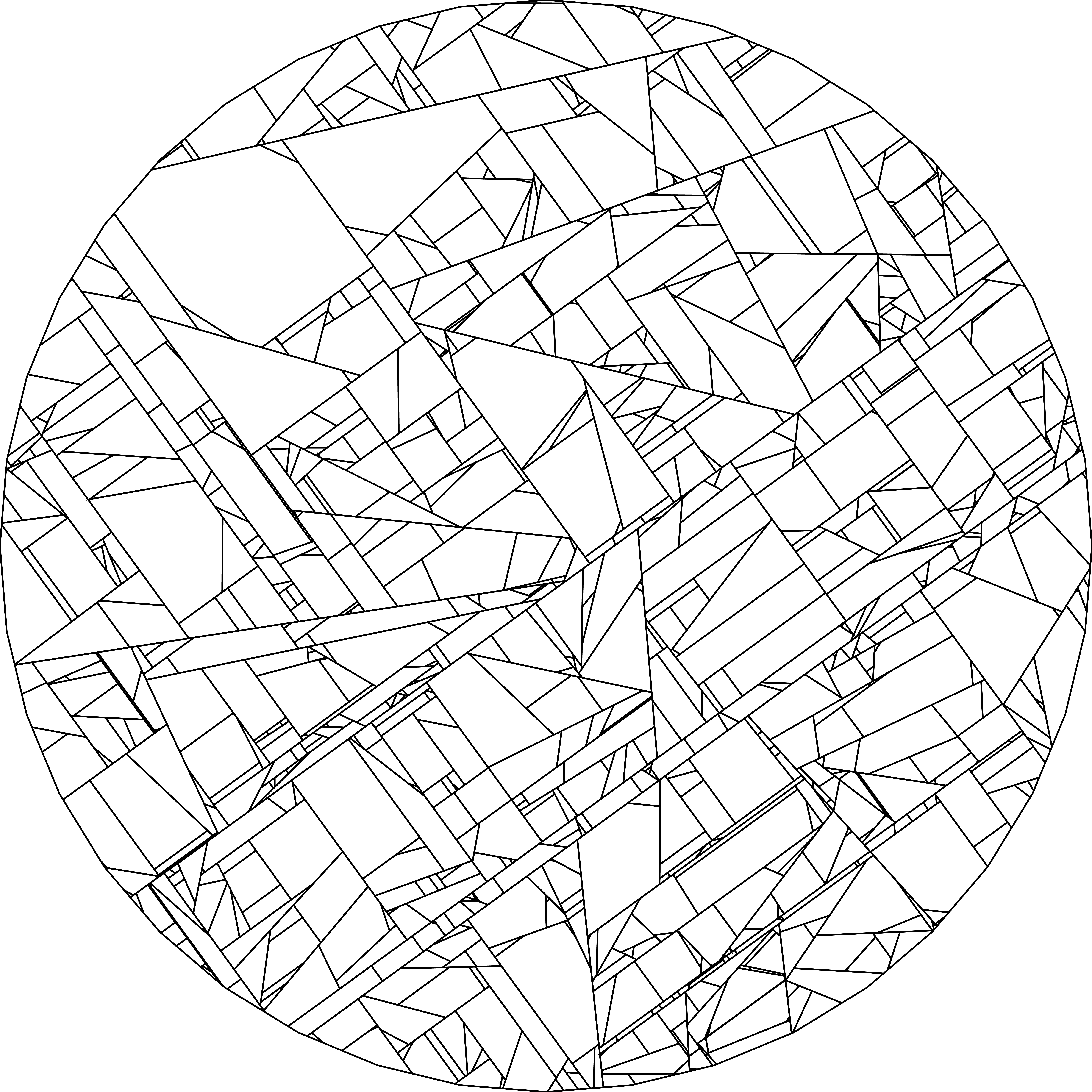}
                
      \end{subfigure}
      \begin{subfigure}[b]{0.3\textwidth}
                \centering
               \includegraphics[width =\textwidth]{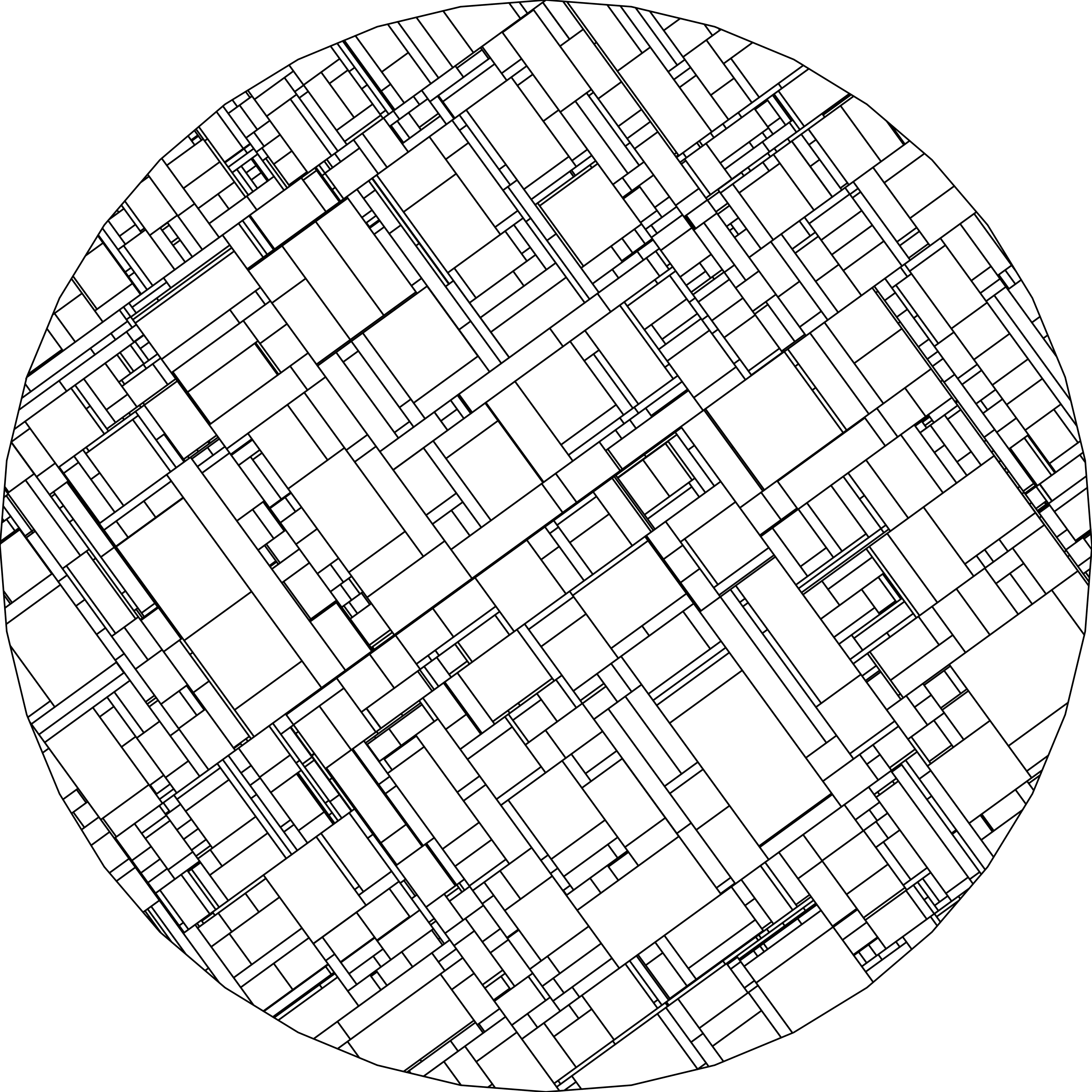}
                
      \end{subfigure}
\end{center}

\caption{Realization of tessellations into a disc. The first row shows PLT with from left to right an anisotropy of $0$, $0.5$ and $1$. The second row shows Crack STIT with the same anisotropy distributions.}
\label{fig:tessellations}
\end{figure}

\paragraph*{Model choice}
Poisson Vorono\"{i} Tessellations as other theoretically well known models have been suggested to represent street systems \cite{Gloaguen2002,Gloaguen2006}.
\\ 
Two models are particularly relevant to model street systems: Poisson Line Tessellation \cite{Stoyan1995,Morlot2012,Kendall2010} and Crack STIT Tessellation \cite{Nagel2005}. 
\\	
It has been shown in \cite{Courtat2011,Perna2011} that a city's morphology can result from two growth mechanisms. The first one is an organic like growth: independent agents divide sequentially plots of land to settle in. Agents do not consult each other, consequently axis they create are not coherent and street intersections are T-shaped. Conversely in the case of a planned growth, agents act under the authority and draw long transportation axis to optimize displacements within the city. 
Crack Stit tessellations mimic the organic case and Poisson Line Tesselations the planned one (Fig.\ref{fig:tessellations}).  
 \\ 
 Moreover these two models are quite similar and can be treated in the same framework. Indeed they both result from sequential divisions of polygons and if edges in a  Crack Stit were extended into infinite line, one would obtain a Poisson Line Tessellation.  
This has for consequence that the typical cell of the tessellations are equal in distribution (\cite{Nagel2008}).

They both depend on two parameters: an intensity parameter $\lambda \in \mathbb{R}_+$ and a probability measure over $[0,\pi]$, $\mathcal{R}$ that describes the anisotropy of the street system.

\paragraph*{Poisson Line Tessellation}
\begin{defin}[Poisson Line Tessellation]
If $\partial \Xi$ is a Poisson Line Process then the Poisson Line Tessellation (PLT) $\Xi$ is the  set of connected components of $\mathbb{R}^2 \backslash \partial \Xi$.
\end{defin}
\paragraph*{Crack STIT Tessellation}
A Crack is the result of a sequential random division of the plane. \\
Informally, at a time $t$, the tessellation is $\mathbf{\Xi}_t = \{ C_i \}$ and  $dt$, each cell $C_i$ in the tessellation has a probability $\lambda.\nu(C_i).dt, \, \lambda>0$ to be divided into two new cells and a probability $o(dt^2)$ to be divided twice. $\nu(.)$ is a positive measure on the set of convex bodies, invariant under rigid motion (for instance area, perimeter, number of vertices)  \cite{Cowan2010}). The tessellation is observed at a finite time $\tau$, the homogeneous quantity $\lambda.\tau$ describes the intensity of the process in the such a way one can come down to $\tau = 1$. If the measure $\nu$ is the perimeter (which is the case under consideration in what follows), the resulting tessellation process has interesting properties: it is STable under ITeration (STIT) and its typical cell is equal in distribution to the PLT's one.   

\begin{defin}[Binary tree]
Let $0,+,-$ be three symbols and let's define recursively the sets $A_n$ by: $A_0 = \{(0)\}$ and $\forall n > 0, \, A_n = \{(a, +)\} \cup \{(a,-)\ , \, a \in A_{n-1}\}$.\\
The rooted binary tree $A$ is defined by $A = \cup A_n$.
\\ 
If $a \in A$ then $a_+ =(a,+) \in A$ and $a_-=(a,-) \in A$ are the daughters of $a$. 
\\
Conversely, if $b \in A, b \neq 0, \exists a \in A, b$ is a daughter of $a$. 
\end{defin}
\begin{defin}[Finite division process]
Let $W$ a compact set of the plane $A$ the rooted binary tree and $\lambda, \tau \in \mathbb{R}_+^*$ 
\\
 $c$ is a random polygon process, $l$ a random line process, $T$ and $\Delta T$ random real number process ; they are all indexed by the elements of $A$ and are defined recursively:

\begin{equation}
\begin{cases}
c(0) = W, T(0)=0 \\
\Delta T(a) \sim \mathcal{E}\left( \lambda \vert c(a) \vert \right), \, \forall a \in A \\ 
l(a) \sim \mathbb{U}_{C(T(a))}, \, \forall a \in A \\ 
c(a_+) = c(a) \cap l(a)_+\, \forall a \in A \\ 
c(a_-) = c(a) \cap l(a)_-\, \forall a \in A \\ 
T(a_+) = T(a_-) = T(a) + \Delta T(a)\, \forall a \in A
\end{cases}
\end{equation}
from these processes, we define the stochastic polygon process $C$ indexed by $A \times \mathbb{R}_+$
\begin{equation}
C(a,t) = 
\begin{cases}
\emptyset & \text{if } t < T(a) \\
c(a) & \text{if } T(a) \leq t < T(a) + \Delta T(a) \\
\emptyset & \text{if } t \geq T(a) + \Delta T(a)
\end{cases}
\end{equation}
The finite division process of $W$ is the tessellation $\Xi_{\tau} = \cup_{a \in A} C(a, \tau)$
\end{defin}

\begin{defin}[Crack STIT Tessellation] \label{prop:crack}
There exists a stationary locally finite tessellation in the plane such as $\forall W \text{ compact subset, }\Xi \cap W \egdist C_{\tau}(W)$ \cite{Nagel2005}. 
\\
It is called the Crack STIT Tessellation. 
\end{defin}

\begin{propi}
Both PLT and Crack STIT Tessellations are locally finite and stationary.
\end{propi}

\subsection{Mean formulae}
Mean formulae (Tab.\ref{tab:formulesMoyennes}, \cite{Nagel2008}) are known for topological features of PLT and Crack STIT in a disc of area $1$ in function of their intensity $\lambda$ and of the so called anisotropy parameter:   
\begin{equation} 
\xi \egdef \iint |\sin{\measuredangle (u,v) }| \mathcal{R}(du)\mathcal{R}(dv) 
\end{equation}
 (these formulae remain true for all Borelian of area $1$ if the tessellation is isotropic i.e. if $\mathcal{R}$ is uniform or if $\xi = 2/\pi$). 
 \\
 Knowing these mean formulae permits to calibrate the models to fit real data which is of prim interest. For instance in sections \ref{sec:propagation} and \ref{sec:resultats} we will chose the intensity of the tessellation such as a block of houses has a perimeter of $400 \, \text{m}$ in average. 
\begin{table}[h!]
\begin{tabular}{cc|c|c}
\hline 
Parameters & Notation & Mean value per u.a for PLT & Mean value per u.a for Crack \\ 
\hline 
Total edge length & $L_A$ & $\lambda$ & $\lambda. \tau $ \\ 

Number of vertices & $N_0$ & $\frac{1}{2} \xi \lambda^2 $ & $L_A^2 \xi$ \\ 
 
Number of edges & $N_1$ & $\lambda^2\xi$ &  $\frac{3}{2}L_A^2 \xi$ \\ 
 
Number of cells & $N_2$ & $\frac{1}{2}\lambda^2\xi$ &  $\frac{1}{2}L_A^2 \xi$ \\ 
 \hline 
Length of the typical edge & $U_1$ & $2/(3\lambda\xi)$ & $ 2/(3L_A\xi)$ \\ 

Perimeter of the typical cell & $U_2$ & $4/(\lambda\xi)$ & $ 4/(L_A\xi)$ \\ 

Area of the typical cell & $A_2$ & $2/(\lambda^2\xi)$ & $ 2/(L_A^2\xi)$ \\ 
\hline 
\end{tabular} 
\caption{\label{tab:formulesMoyennes} Expectancies of various morphological features of PLT and Crack STIT in function of their intensity $\lambda$ and their anisotropy parameter $\xi$. The notion of "typical" object is very useful in stochastic geometry. It basically corresponds to the idea of an object sampled out of a large collection independently of any measure that can be applied to it (see \cite{Baccelli2009,Stoyan1995,Lieshout2000}).  }
\end{table}

In a city modelling context, the family of angular distributions:

\begin{equation} \mathcal{R}_{\rho} = \rho\mathcal{U}_{[0,\pi]}+ \frac{(1-\rho)}{2}(\delta_{\theta} + \delta_{\theta + \pi/2}) 
\end{equation}.
allows to go continuously from an isotropic network ($\rho = 0$) to an anisotropic Manhattan-like one ($\rho = 1$). For this family of distributions, $\xi$ writes: 

\begin{equation}
\xi(\rho) = \rho^2 \left(\frac{1}{2} - \frac{2}{\pi}\right)+ \frac{2}{\pi} 
\end{equation}  

\section{Building generation \label{sec:batiments}}

\begin{figure}[h!]
\begin{center}
\begin{subfigure}[b]{0.3\textwidth}
                \centering
                \includegraphics[width =\textwidth]{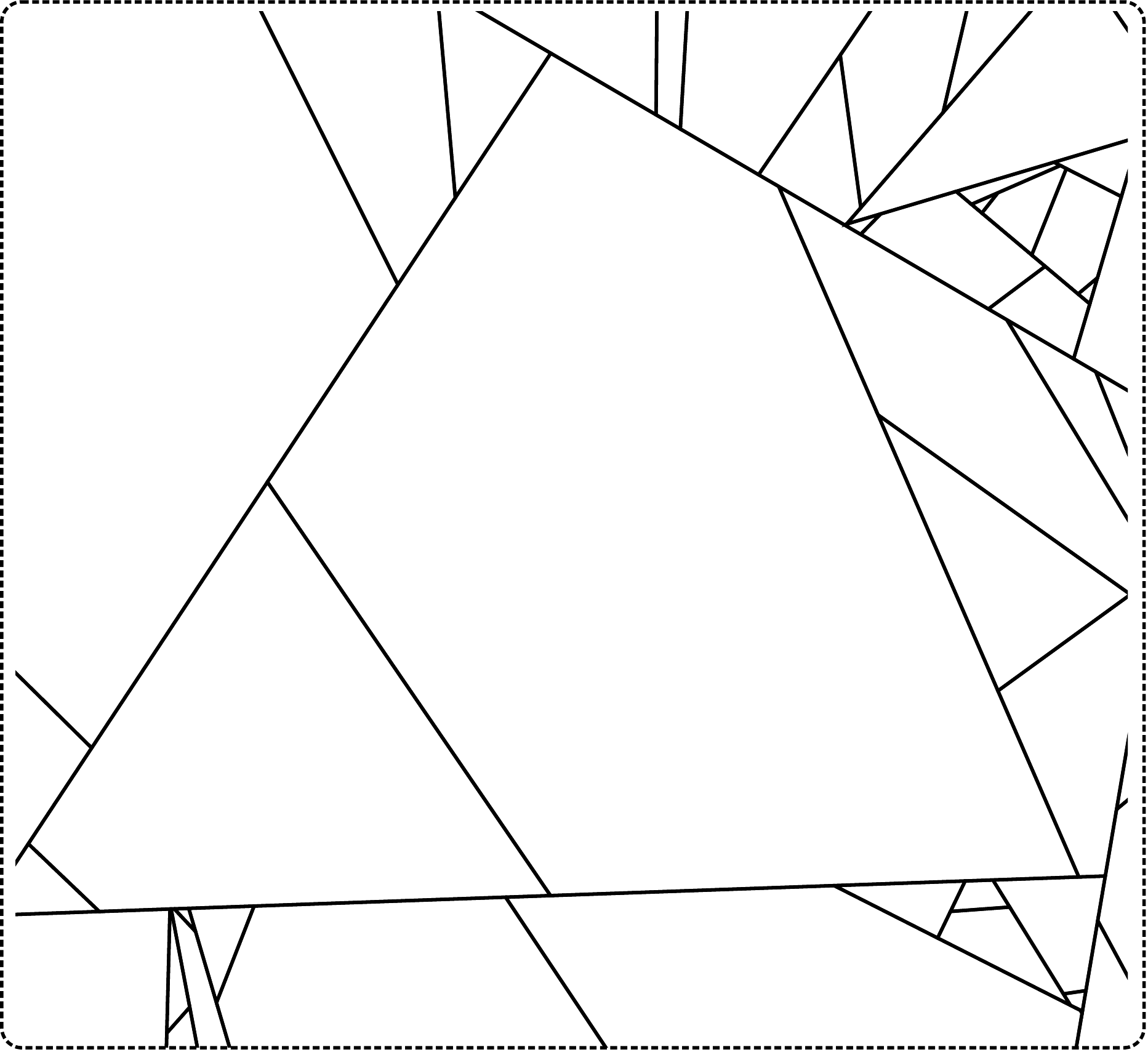}
              
      \end{subfigure}
      \begin{subfigure}[b]{0.3\textwidth}
                \centering
                \includegraphics[width =\textwidth]{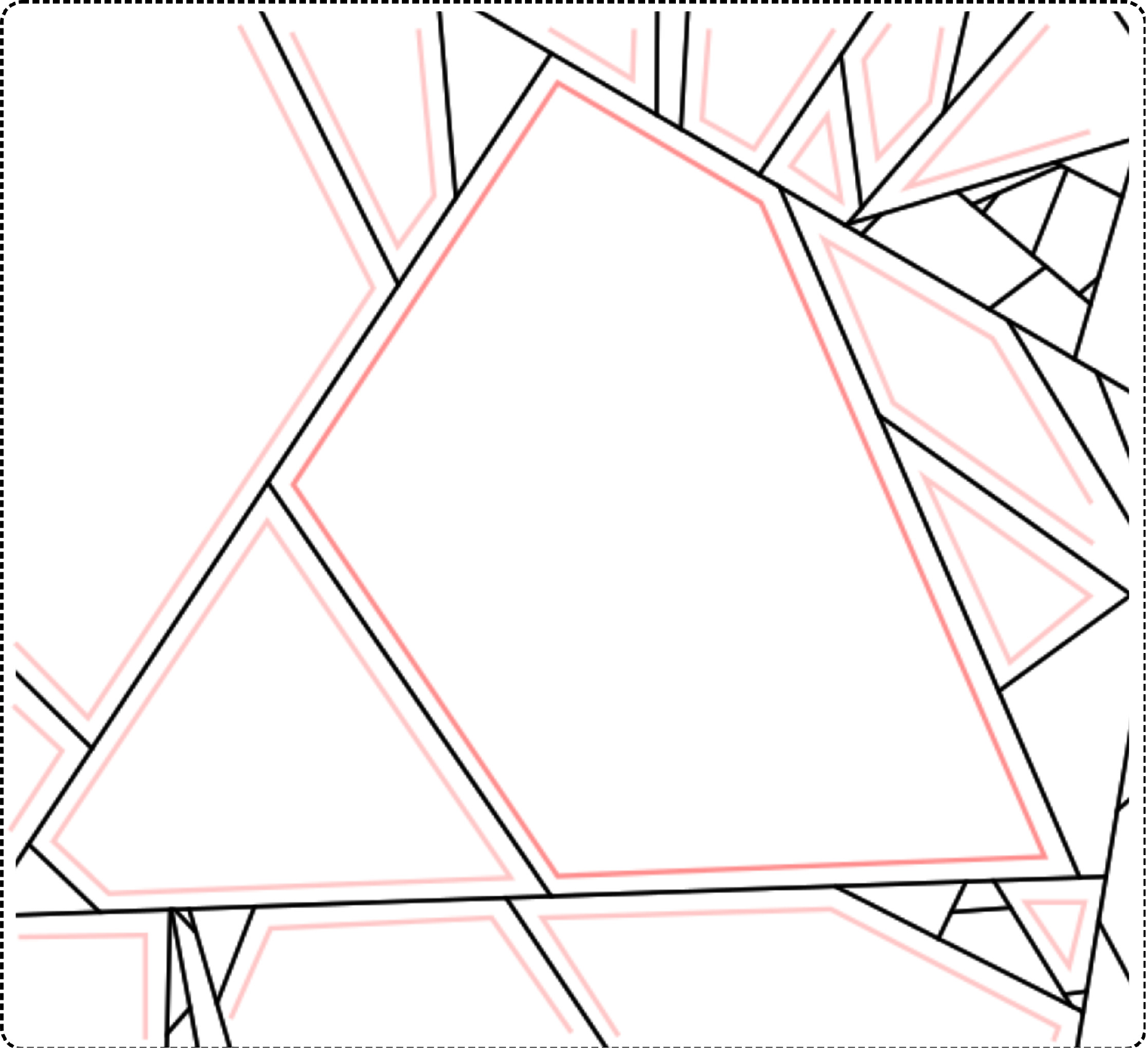}
               
      \end{subfigure}
      \begin{subfigure}[b]{0.3\textwidth}
                \centering
                \includegraphics[width =\textwidth]{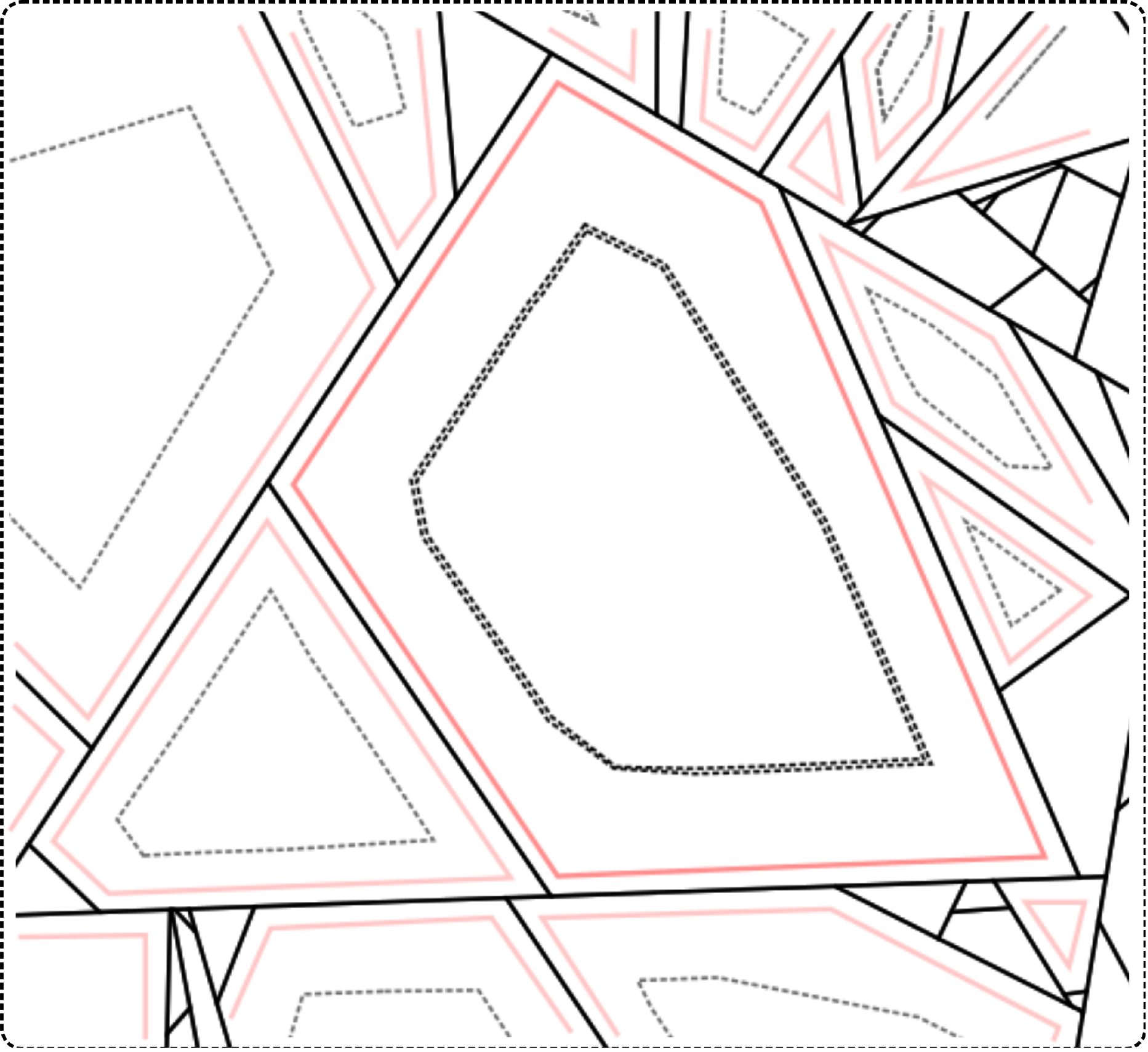}
                
      \end{subfigure}
      \begin{subfigure}[b]{0.3\textwidth}
                \centering
                \includegraphics[width =\textwidth]{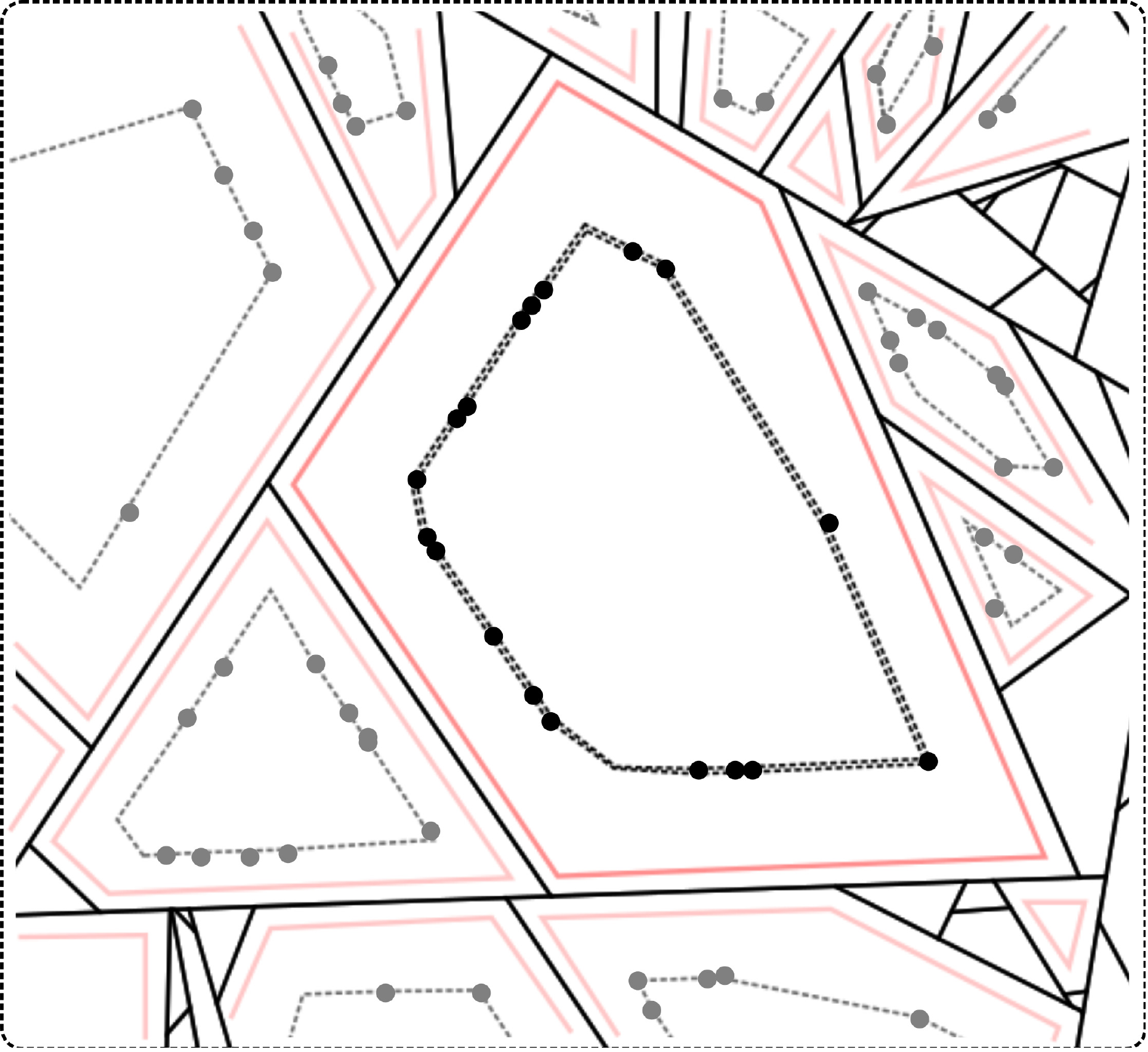}
              
      \end{subfigure}
      \begin{subfigure}[b]{0.3\textwidth}
                \centering
                \includegraphics[width =\textwidth]{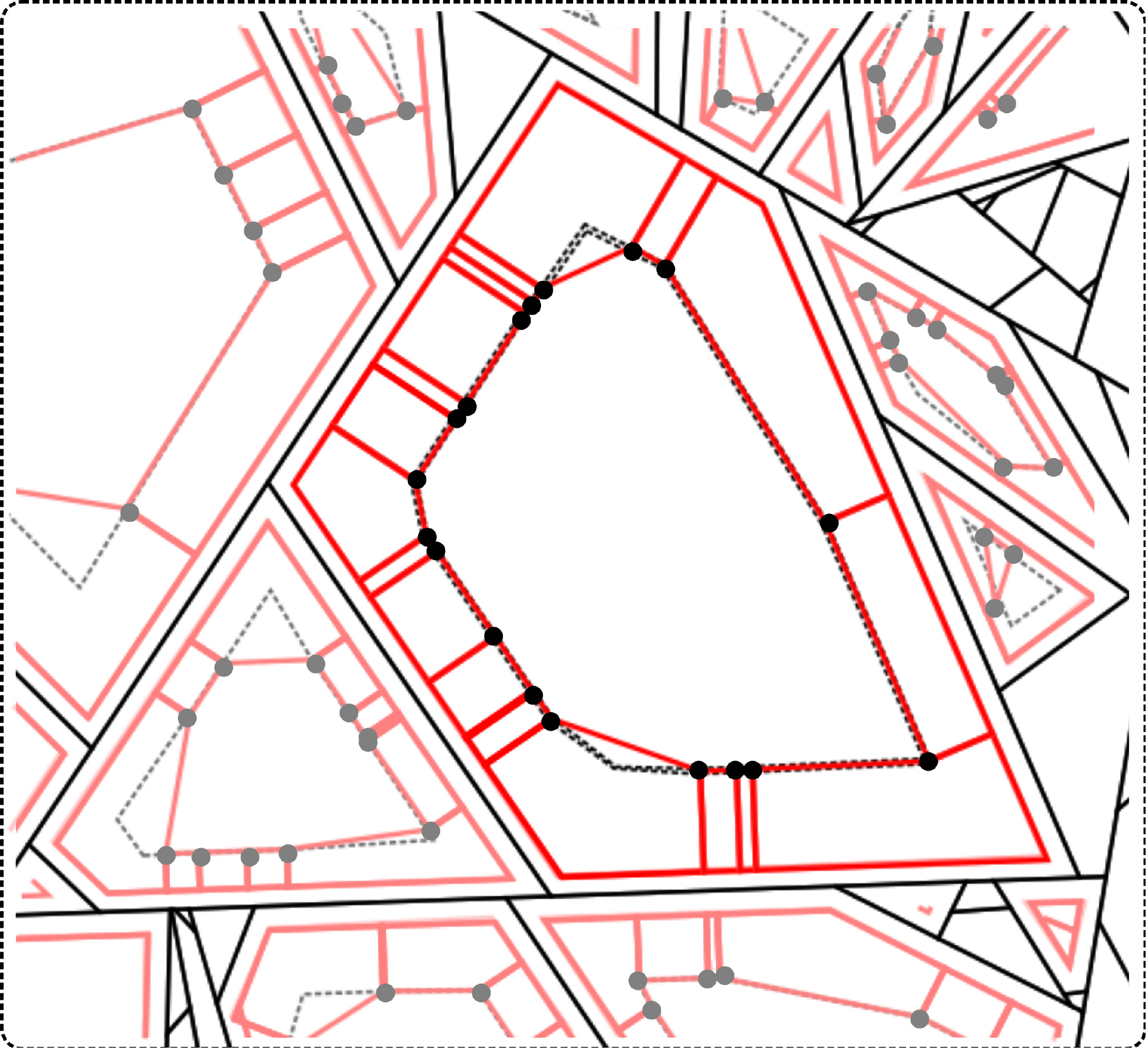}
                
      \end{subfigure}
      \begin{subfigure}[b]{0.3\textwidth}
                \centering
                \includegraphics[width =\textwidth]{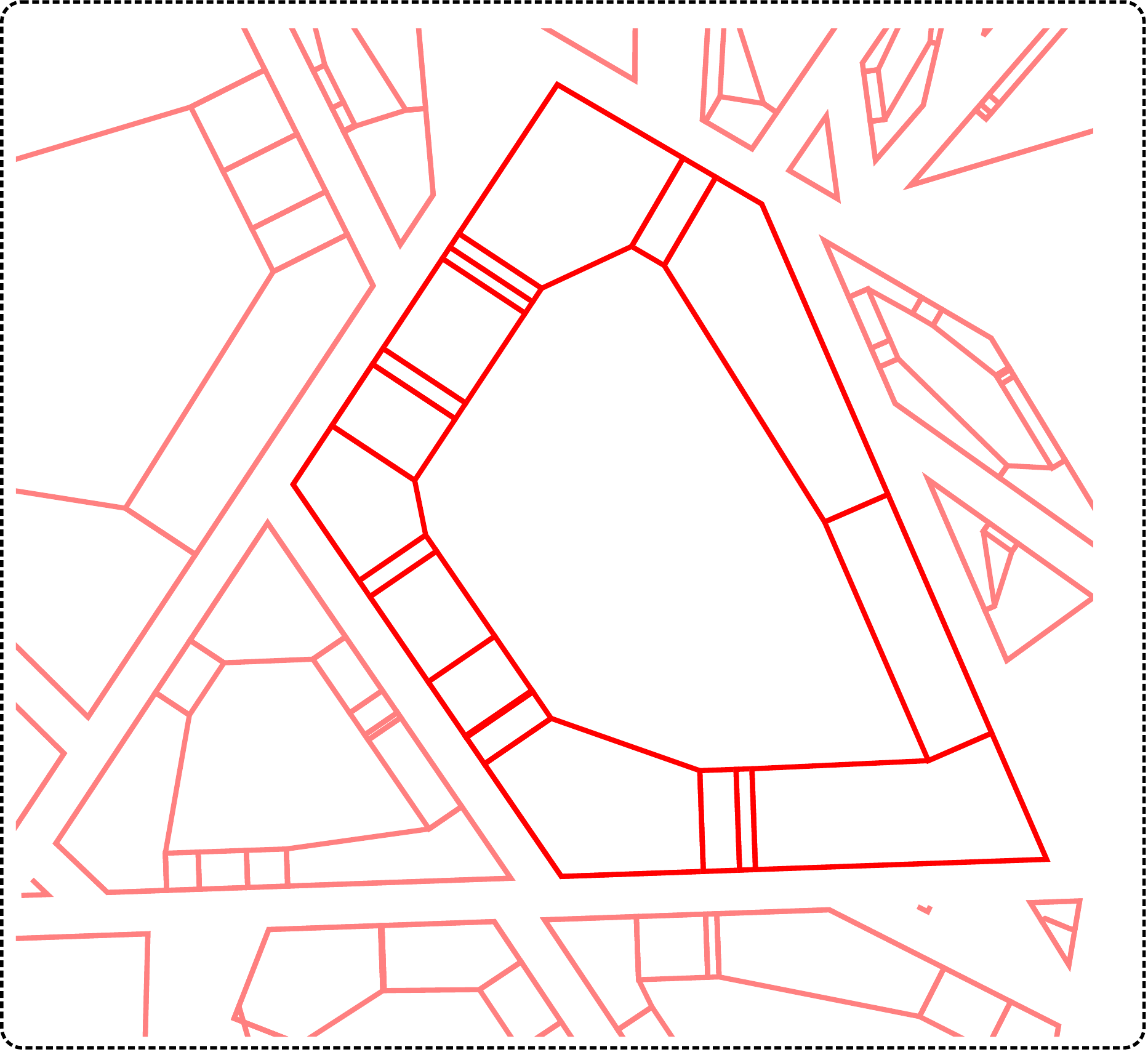}
               
      \end{subfigure}
\end{center}

\caption{\label{fig:construction} Steps in the building generation. From a tessellation (1) we apply an erosion operator to axis (2) in each new cell, we compute its dilated polygon with respect to its center of mass (3) we draw on this polygon a Poisson Point Process (4) whose points are projected to create buildings' footprint (5)  }
\end{figure}

The tessellation is to represent the street axis (alignments of edges), the city skeleton (Fig.\ref{fig:construction},1). From this, axis are thickened with a Minkowski's sum $\oplus_{\epsilon}, \, \epsilon >0 $. 
\begin{defin}
If $A$ is a subset of the plane, 
\begin{equation}
 \oplus_{\epsilon} A = \{ x \in \mathbb{R}^2, d(x,A) \leq \epsilon \}
\end{equation}
\end{defin}

\begin{propi}
If $\{ A_i\}$ are the axis of a tessellation, the connected components of $\mathbb{R}^2 \setminus \oplus_{\epsilon} \cup A_i, \, \epsilon >0 $ are polygons $B_i$ that do not intersect.
\\
If $\{ C_i\}$ is the set of cells of the tessellation then each $B_i$ is the image of a cell $C_k$ by the operator 
\begin{equation}
 \ominus_{\epsilon / 2}(C) = \{ x \in C ,  d(x, \partial C)  >	 \frac{\epsilon}{2}\}
 \end{equation}
\end{propi}

As a consequence, thinkening axis into streets is equivalent to generate building blocks from cells (Fig.\ref{fig:construction},2) by applying  $\ominus_{\epsilon / 2}$ independently to each cell in $\mathcal{T}$.

Once blocks obtained, we associate to each block $B$ its image by the dilatation of center its center of mass and ratio $\eta$: $\tilde{B}$ (Fig.\ref{fig:construction},3). 

A Poisson Point Process of intensity  $1/ (b.\eta) $ is drawn on $\partial \tilde{B}$ (Fig.\ref{fig:construction},4).Then we project theses points on the border of $B$ and apply a computational procedure to  create polygons from these two sets of points. Simple but laborious to explain symbolically, this procedure is shown in (Fig.\ref{fig:construction},5) and permits to draw buildings' footprint.
\\
To each building is associated a random height from a distribution $ \mathcal{E}(h), \, h>0$  
 
From these tessellations streets models and this building generation method, it is possible to simulate a random urban environment by deciding of:
\begin{itemize}
\item It's topology (mainly T intersections of X intersections) 
\item A mean parameter chosen from Tab.\ref{tab:formulesMoyennes} to fit the intensity of the model,
\item The mean street width,
\item The mean number of buildings in a block,
\item The mean building height.
\end{itemize}

These random models define parametric classes of urban environments morphologicaly equivalent. We present in \ref{ann:tess} key functions to simulate efficiently a large number of representative environments for each class and set of parameters. 
\\ 
In what follows, we illustrate the interest of these models by studying the environment's mean morphology impact over the path loss function of an electromagnetic wave emanating from a placed-high antenna to any user at the ground level, at a distance $d$ of the antenna. The commonly admitted rule of thumb is that $P(d) \simeq d^{-\alpha}$ with $\alpha 0$ depending on the environment.
\\ 
In \ref{sec:propagation} we settle a probabilistic framework suited to compute the wave propagation by Monte-Carlo / ray tracing methods. 
\\ 
In \ref{sec:resultats} we implement these methods for several representative morphology classes and a large number of independently generated environments and exhibit statistically the impact of the morphology on the path loss exponent $\alpha$.

Given a real world map, one can fast assesses fit the map to a statistical model and deduce the corresponding expected $\alpha$.
\section{Propagation simulation  \label{sec:propagation}}

\subsection{Source antenna model}
Given a urban environment, we model a wireless telecommunication antenna by a sphere $\mathbb{S}$, placed at the center off mass of the nearest roof of the origin of the plane ($\vec{0}$), with an additional height $\delta H$.

\begin{figure}[h!]
\begin{center}

                \includegraphics[width =\textwidth]{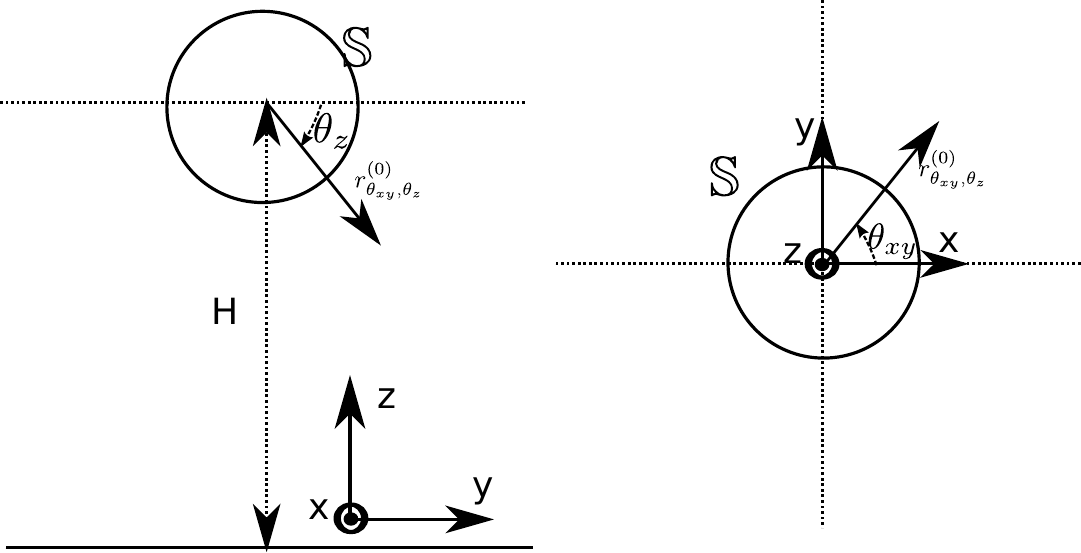}

      \caption{\label{fig:source} Sketch of the source $\mathbb{S}$ as a sphere observed from the front and from above. One can see the definition of angles $\thz \, \thx$ that identify the first portion of a ray $r^{(0)}_{\thx, \thz}$.  }
      \end{center}
      \end{figure}

It is parametrized with spherical coordinates $s=(\thx, \thz)$ (Fig.\ref{fig:source}), equipped with its Borel $\sigma$-algebra $\tribS$ and the uniform measure $ \mathcal{S}(\dif \thx, \dif \thz) = \dfrac{1}{4\pi} \dif \thx|\sin{(\dif \thz)}|$.
\\
It emits in the sphere portion $\Delta\mathbb{S} = [\thx_0-\delta \thx, \thx_0+\delta \thx ] \times [\thz_0-\delta \thz, \thz_0+\delta \thz ]$ of total measure $\mathcal{S}(\Delta \mathbb{S}) = \dfrac{1}{\pi}\delta \thx  \sin{(\delta \thz)}\cos{(\thz_0)} $. The measure restricted to $\Delta \mathbb{S}$ is $\Delta \mathcal{S}(.) \egdef \dfrac{\mathcal{S}(.)}{\mathcal{S}(\Delta \mathbb{S})}$. The power density of the antenna is $P_0.\Delta \mathcal{S}(.)$ with $P_0 \geq 0$ the total power emanating from the source. 
      
We write $\partial E$ the set of polygons constituting the buildings and the ground of the city. $\partial E$ splits the space $\Rt$ into three parts: $E$ which is the interior of the buildings and the under-ground where wave cannot propagate, $\partial E$ itself where waves are reflected and $\Omega = \Rt \backslash (E \cup \partial E)$, an open set whose boundary is $\partial E$, representing the wave propagation space. 
\\
To each $s \in \Delta \mathbb{S}$ is associated $r_s \subset \Omega$ the trajectory of an elementary wave emanating from $s$. It is assumed that $r_s$ is the result of successive reflections on $\partial E$ respecting the Snell-Descartes rules. 
\\
Between its $n$-th and $n+1$-th reflection, $r_s$ is a segment denoted $r_s^{(n)}$ ; $r_s =\cup_n r_s^{(n)}$. If $x \in \Omega \cap r_s^{(n)}$, we write $\vec{r}_s^{(n)}(x)$ the propagation direction of $r_s$ at point $x$. 
\\
When reflected, a portion $1-\gamma, \, \gamma \in [0,1]$ of the ray power is absorbed.
 
We define in the following subsection in probabilistic terms the physical notions of electromagnetic flow, Poynting vector and electromagnetic power.  

\subsection{Flow, Poyning's vector and power probabilistic definition}

We emphasize that $\ensS$ is compact, separable and connected metric space. 
\\
The following assumption is a convenient abstract description of the trajectories of rays from the source $\ensS$: 
\begin{hyp} \label{hyp:1}
Let $n\in \N$, $x_0 \in \Omega$ and $s_0 \in \ensS$. 
\\ 
There exist $\epsilon_{x_0}>0$ and $\epsilon_{s_0} >0$ and balls $V_{\epsilon_{x_0}}(x_0) \subset \Omega$, centered in $x_0$, of radius $\epsilon_{x_0}$ and $V_{\epsilon_{s_0}}(s_0) \subset \ensS$, centered in $s_0$, of radius $\epsilon_{s_0}$ such as 
\begin{itemize}
\item If $x_0 \in r_{s_0}^{(n)}$, $\left\lbrace r_s^{(n)} \cap V_{\epsilon_{x_0}}(x_0) , \, s \in V_{\epsilon_{s_0}}(s_0)   \right\rbrace$ is homeomorph to $V_{\epsilon_{s_0}}(s_0) \times ]0,1[$
\item If $x_0 \notin r_{s_0}^{(n)}$, $\left\lbrace r_s^{(n)} \cap V_{\epsilon_{x_0}}(x_0) , \, s \in V_{\epsilon_{s_0}}(s_0)   \right\rbrace = \emptyset$
\end{itemize}
\end{hyp} 
\begin{propi}
If $n\in \N$, $x_0 \in \Omega$ then the set $\left\lbrace s \in \ensS , \, x\in r_s^{(n)}\right\rbrace$ is finite.
\end{propi}
In what follows, without more precision, $\ensM$ is a separable, connected, locally compact manifold in $\Omega$ of dimension $\dim \ensM \leq 2$, $C^1$ almost everywhere, equipped with its Borel $\sigma$-algebra $\tribM$ and its Borel measure $\mesM$. 

\begin{defin} \label{def:ens1}
For all $n \in \N$, $m \in \N \cup \{ \infty \}$,  $s \in \ensS$ and $(S,M) \in \tribS \otimes \tribM$, we define the sets: 
\begin{itemize}
\item $q_{\tribM}^{(n)}(s) = \left\lbrace x \in \tribM , x \in r_s^{(n)} \right\rbrace$
\item $ R_{\ensM}^{(n)}(S,M) = \left\lbrace s \in S, q_{\tribM}^{(n)}(s) \neq \emptyset \right\rbrace$
\item $ R_{\ensM}^{(n,m)}(S,M) = \left\lbrace s \in S, \sharp q_{\tribM}^{(n)}(s) =m \right\rbrace$
\end{itemize}
\end{defin}

\begin{propi}\label{propi:1}
$\forall n \in \mathbb{N}, \, \mathcal{S}\left(R_{\mathbb{M}}^{(n,\infty)}\right)=0$
\end{propi}
\begin{propi}\label{propi:2}
If $\dim \ensM \leq 1$, $\forall n \in \mathbb{N}, \, \forall m \in \mathbb{N}, \, \mathcal{S}\left(R_{\mathbb{M}}^{(n,m)}\right)=0$
\end{propi}
\begin{defin}
We define for $(S,M) \in \tribS \otimes \tribM$,  
\begin{equation}
\left\vert \Phi_{\mathbb{M}}^{(n)} \right\vert(S,M) = \sum_m m\cdot\mesS \left( R_{\mathbb{M}}^{(n,m)}(S,M) \right)
\end{equation}
\begin{itemize}
\item $\left\vert \Phi_{\mathbb{M}}^{(n)} \right\vert$ is a measure over $\ensS \times \ensM$
\item Moreover, $ \left\vert \Phi_{\mathbb{M}}^{(n)} \right\vert \ll \mesS \otimes \mesM$
\end{itemize}
In particular, if $S$ is an open set in $\ensS$, $\left\vert \Phi_{\mathbb{M}}^{(n)} \right\vert(S,\cdot)$ is a measure over $ \mathbb{M}$ and $\left\vert \Phi_{\mathbb{M}}^{(n)} \right\vert(S,\cdot) \ll \mesM$ 
\end{defin}
\begin{defin}
Let's write $\vec{m}(x)$ the orthogonal direction of $\ensM \ni x$ at point $x$.
\\
In addition to \ref{def:ens1}, we define for all $n \in \N$, $m \in \N \cup \{ \infty \}$,  $s \in \ensS$ and $(S,M) \in \tribS \otimes \tribM$, the sets: 
\begin{itemize}
\item $q_{\tribM, +}^{(n)}(s)= \left\lbrace x \in r_s^{(n)} \cap \tribM, \left\langle \vec{r}_s^{(n)}(x) \middle\vert \vec{m}(x) \right\rangle >0   \right\rbrace$
\item $q_{\tribM, -}^{(n)}(s)= \left\lbrace x \in r_s^{(n)} \cap \tribM, \left\langle \vec{r}_s^{(n)}(x) \middle\vert \vec{m}(x) \right\rangle <0   \right\rbrace$
\item  $ R_{\ensM,+}^{(n,m)}(S,M) = \left\lbrace s \in S, \sharp q_{\tribM,+}^{(n)}(s) =m \right\rbrace$
\item  $ R_{\ensM,-}^{(n,m)}(S,M) = \left\lbrace s \in S, \sharp q_{\tribM,-}^{(n)}(s) =m \right\rbrace$
\end{itemize}
\end{defin}
\begin{defin}
We define for $(S,M) \in \tribS \otimes \tribM$,  
\begin{itemize}
\item $ \Phi_{\mathbb{M}+}^{(n)} (S,M) = \sum_m m\cdot\mesS \left( R_{\mathbb{M}+}^{(n,m)}(S,M) \right)$
\item $ \Phi_{\mathbb{M}-}^{(n)} (S,M) = \sum_m m\cdot\mesS \left( R_{\mathbb{M}-}^{(n,m)}(S,M) \right)$
\end{itemize}
and 
\begin{equation}
\Phi_{\mathbb{M}}^{(n)} = \Phi_{\mathbb{M}+}^{(n)} - \Phi_{\mathbb{M}-}^{(n)}
\end{equation}
\end{defin}
\begin{propi}
$\Phi_{\mathbb{M}}^{(n)}$ is a signed measure over $\mathbb{S}\times\mathbb{M}$.
\\
Its total variation is $\vert \Phi_{\mathbb{M}}^{(n)} \vert$ and $\forall n \in \mathbb{N}, \quad \Phi_{\mathbb{M}}^{(n)} \ll \mathcal{S} \otimes \mathcal{M}$. 
\\ 
We write $\varphi_{\mathbb{M}}^{(n)} = \frac{d\Phi_{\mathbb{M}}^{(n)}}{d\mathcal{S} \otimes \mathcal{M}}$.
\end{propi}
\begin{rem}
If $S$ is an open set in $\ensS$, $\Phi_{\mathbb{M}}^{(n)}(S, \cdot)$ is a signed measure over $\mathbb{M}$.
\\
Its total variation is $\vert \Phi_{\mathbb{M}}^{(n)}(S, \cdot) \vert$ and $\forall n \in \mathbb{N}, \quad \Phi_{\mathbb{M}}^{(n)}(S, \cdot) \ll  \mathcal{M}$.
\\ 
We write $\varphi_{\mathbb{M}}^{(n)}(S, \cdot) = \frac{d\Phi_{\mathbb{M}}^{(n)}(S, \cdot)}{d\mathcal{M}}$
\end{rem}
\begin{rem}
$\varphi_{\mathbb{M}}^{(n)}$ is degenerated:
\\
if $x \in \Omega$, $\left\lbrace s \in \ensS, x \in r_s^{(n)} \right\rbrace $ is finite. By writing $K(x)$ the size of this set and $s_1(x), \cdots ,s_{K(x)}(x)$ its elements, there exist $S_1(x), \cdots , S_{K(x)}(x)$ open sets of $\ensS $ such as: 
\begin{equation}
\varphi_{\mathbb{M}}^{(n)}(s,x) = \sum_{k = 1}^{K(x)} \varphi_{\mathbb{M}}^{(n)}(S_k(x), x) \delta(s - s_k)
\end{equation}
\end{rem}
\begin{propi}
If $S$ is an open set in $S$ and if $V_{\epsilon}(x)$ is a neighborhood system of $x$ indexed by $\epsilon >0$, with $\lim_{\epsilon \to 0} \text{diam}(V_{\epsilon}(x))  =0$ then 
\begin{equation}
 \varphi_{\mathbb{M}}^{(n)}(S,x) = \lim_{\epsilon \to 0} \frac{\Phi_{\mathbb{M}}^{(n)}(S, V_{\epsilon}(x))}{ \mesM( V_{\epsilon}(x))}
 \end{equation}
\end{propi}
\begin{propi}
If $\ensM$ is a closed surface and $S$ an open set in $\ensS$, $\forall n \in \N, \, \Phi_{\mathbb{M}}^{(n)}(S, \ensM) =0$.
\end{propi}
\begin{propi}
Let $S$ an open set in $\ensS$.
\\
There exist $\varphi_S^{(n)}: \Omega \to \R_+$ and $\vec{d}_S^{(n)}: \Omega \to \ensU^2$ the two dimensional unitary sphere in $\R^3$, continuous function such as forall $x \in \Omega$ and for all affine plane $\ensM \ni x$ orthogonal to the direction $\vec{m} \in \ensU^2$, $\varphi_S^{(n)}$ and $\vec{d}_S^{(n)}$ are integrable over $\ensM$ and 
\begin{equation}
\varphi_{\mathbb{M}}^{(n)}(S,x) = \left\langle  \varphi_S^{(n)}(x) \vec{d}_S^{(n)}(x) \middle\vert \vec{m}\right\rangle
\end{equation}  
If $\varphi_{\mathbb{M}}^{(n)}(S,x) \neq 0$, $\varphi_S^{(n)}(x)$ and $\vec{d}_S^{(n)}(x)$ are unique.
\\ 
Moreover in that case, there is a positive finite number $K \geq 1$ of elements of $S$, $s_1, \cdots, s_K$, for all $k \leq K$, $x \in r_{s_k}^{(n)}$,  each one being associated to an open set $S_k \ni s_k$ included in $S$, such as 
\begin{equation}
\varphi_S^{(n)} \vec{d}_S^{(n)} = \sum_{k=1}^K \varphi_{S_k}^{(n)}(x)\vec{r}_{s_k}^{(n)}(x)
\end{equation}
\end{propi}
\begin{defin} \label{def:global}
Let $\gamma \in \R, \vert \gamma \vert <1$. 
\\
For all $x \in \Omega$, we define: 
\begin{itemize}
\item The total flow: $\Phi_{\ensM}(M) = \sum_n \gamma^n \Phi_{\ensM}^{(n)}(\Delta \ensS, M)$
\item The elementary flow: $\varphi_{\ensM}(x) = \sum_n \gamma^n \int_{\Delta \ensS} \varphi_{\ensM}^{(n)}(s,x) \mesS (\dif s)$
\item The Poyning vector: $\vec{\Pi}(x) = \sum_n \gamma^n \varphi_{\ensS}^{(n)}(x) \vec{d}_{\ensS}^{(n)}(x)$
\item The power: $P(x) = \left \Vert \vec{\Pi}(x) \right\Vert$
\end{itemize}
\end{defin}
\begin{propi}
The quantities defined in \ref{def:global} are bonded by the relationships: 
\begin{itemize}
\item $\Phi_{\ensM} \ll \mesM$ and $\dfrac{\dif \Phi_{\ensM}}{\dif \mesM} = \varphi_{\ensM}(x) $
\item If $\vec{m}(x)$ is the orthogonal direction to $\ensM$ at point $x$, $\varphi_{\ensM}(x) = \left\langle \vec{\Pi}(x) \middle \vert \vec{m}(x)  \right\rangle$
\end{itemize}

\end{propi}

\subsection{Power estimation}
\begin{defin}[Power estimator]
Let $x \in \Omega$ and $M \ni x$ an affine plane orthogonal to some direction $\vec{m}$ and $V_{\epsilon}(x)$ a neighborhood system of $x$ in $\ensM$ with $\text{diam} V_{\epsilon}(x) \to_{\epsilon \to 0} 0$.
\\
Let $N \in \N^*$ and $s_1,\cdots, s_N \in \Delta \ensS$, $N$ independent random variables sampled from  $\Delta \mesS$.
\\
$\forall i \leq N$, we write $q^{(n)}(i) =\{ \vec{d}^{(n)}_i\} $ if $r_{s_i}^{(n)} \cap V_{\epsilon}(x) \neq \emptyset $, $q^{(n)}(s_i) = \emptyset$ otherwise.
\\ 
The power estimator in $x$, $\hat{P}_{\epsilon, \vec{m},N}(x)$ is defined by:
\begin{equation}
\hat{P}_{\epsilon, \vec{m}, N}(x) = \dfrac{1}{N\cdot\mesM (V_{\epsilon}(x))} \left\Vert \sum_{i=1}^N \sum_{n=0}^{\infty} \sum_{\vec{d}^{(n)}_i \in q_i^{(n)}} \gamma^n\dfrac{\vec{d}^{(n)}_i}{\left\langle \vec{d}^{(n)}_i \middle\vert \vec{m} \right\rangle} \right\Vert
\end{equation}
\end{defin}
\begin{propi}[Monte-Carlo ray-tracing] 
For almost every direction $\vec{m} \in \ensU^2$, $\hat{P}_{\epsilon,\vec{m},N}$ converges almost surely to $\hat{P}(x)$ as $N\to \infty$ and $\epsilon \to 0$.
\end{propi}
Of course, in practice, the limit is not reachable. We will set $\epsilon$ small compared to the typical variation distance of the problem and $N$ large enough to ensure the resulting variability of the estimator is consistent with the target application.
\\
This variability is investigated in the next section.  
\subsection{Variance reduction}
To assess the accuracy of the previously defined estimators, we place in the simplest configuration where:
\begin{itemize}
\item $\partial E$ is reduced to the ground i.e. the plan of equation $z=0$ also denoted $\ensM$,
\item the antenna $\ensS$ is place above the origin of $\ensM$, at a height $H>0$, 
\item it emits uniformly in a portion $\Delta \ensS$ delimited by $\thx \in [0,2\pi], \, \thz \in ]\thz_0 - \delta\thz, \thz_0+\delta\thz[$, with $\thz_0-\delta\thz>0$ and $\thz_0+\delta\thz \leq \pi/2$, $\delta\thz >0$.
\end{itemize} 
The uniform measure over $\Delta\ensS$ is $\Delta\mesS(\dif \thx, \dif \thz) = \dfrac{1}{C_{\Delta\ensS}}\dif\thx \vert \sin{\dif\thz} \vert$ with $C_{\Delta\ensS} = 4\pi\sin{\delta\thz}\cos{\thz_0}$.

We evaluate the estimators at a ground level point $x \in \ensM$.
\\
For all $s \in \Delta \ensS,$ $r_s$ reaches $\ensM$ one and only one time. Thus the application $r: \Delta \ensS \to \ensM; \quad s \to r_s\cap \ensM$ is an injection.
\\
The image of $\Delta \ensS$ is an open subset of $\ensM$ denoted $\Delta \ensM$ corresponding to the crown centered in $0$, of internal radius $H/\tan{\thz_0+\delta\thz}$ and external radius $H/\tan{\thz_0+\delta\thz}$.
\\ 
For numerical applications, we will consider in what follows $H = 20\text{m}$, the estimations at position $x$ are done from a square $V(x)$ with area $100\text{m}^2$  and $\thz_0, \delta \thz$ such as $\Delta \ensM$ has an internal radius of $50$m and an external radius of $1$km.
\begin{propi}
$r: \Delta \ensS \to \Delta \ensM$ is a $C^1$-differomorphism.
\\ 
It's Jacobian matrix at $s \in \Delta \ensS$ is denoted $J_r(s)$ and $\vert \det{J_r(s=(\thx,\thz))}\vert =H^2\cdot\dfrac{\cos{\thz}}{\sin^3{\thz}}$
\end{propi}
Let $V(x)$ be a (small) Borelian set surrounding $x \in \Delta \ensM$. The flow through $V(x)$ is $\Phi(V(x)) = \Delta \ensS (r^{-1}(V(x)))$ and the power in $x$ is approximate by $P(V(x)) = \dfrac{\Phi(V(x))}{\mesM(V(x))\sin{\thz(x))}}$
\begin{propi}
If $N \in \N^*, \, s_1, \cdots, s_N$ are independent random variable drawn on $\Delta \ensS$ with respect to $\Delta \mesS$ then 
\begin{equation}
\hat{\Phi}(V(x)) = \dfrac{1}{N}\sum_{i=1}^N \mathbf{1}_{V(x)}(r(s_i))
\end{equation}
is an unbiased estimator for $\Phi(V(x))$.
\end{propi}
We will use $\hat{P}(x) \simeq \hat{P}(V(x)) = \dfrac{\hat{\Phi}(V(x))}{\mesM(V(x))\sin{\thz(x))}}$ as an estimator of $P$.
\begin{defin}
If $X$ is a random variable over some probability space with finite expectancy and variance, the relative variation of $X$ is defined by 
\begin{equation}
\sigma_r(X) = \dfrac{\sqrt{\mathbb{V}X}}{\mathbb{E}X}
\end{equation}
\end{defin}
\begin{propi}
The variance of $\hat{\Phi}(V(x))$ is $\mathbb{V}(\hat{\Phi}(V(x))) = \dfrac{\Phi(1-\Phi)}{N}$ ; 
the relative variation of $\hat{\Phi}(V(x))$ and $\hat{P}(V(x))$ are the same and is worth 
\begin{equation}
\sigma_r = \sqrt{\dfrac{1-\Phi}{N\Phi}}
\end{equation}
\end{propi}

In the limit case where $\Phi \ll 1$ and $\mesM(V(x)) \ll H^2$, $\sigma_r \simeq \dfrac{1}{\sqrt{N\Phi}}$ and 
\begin{eqnarray*}
\Phi(V(x)) &\simeq& \cos{r^{-1}(x)}\vert \det{J_r(r^{-1}(x)}\vert^{-1} \dfrac{\mesM(V(x))}{C_{\Delta \ensS}}  \\ 
 \sigma_r &\simeq& \sqrt{\dfrac{C_{\Delta \ensS} H^2}{N\sin^3{x}\mesM(V(x))}} 
 \simeq \Vert x \Vert^{1.5} \sqrt{\dfrac{C_{\Delta \ensS}}{NH\mesM(V(x))}}
\end{eqnarray*}
Numerically, we can afford to simulate up to $N=10^7$ rays this result implies that next to the internal radius of the crown, $\sigma_r \gtrsim 0.4\% $ and next to the external radius, $\sigma_r \gtrsim 33\%$. 
\\ 
In what follows, we design an importance sampling to homogenize the relative variation through space.
\begin{propi}
\label{th:var}
If $\mesT$ is a measure over $\Delta \ensS$ such as $\mesT(\Delta \ensS) = 1 $ and $ \mesS \ll \mesT$, we write $\omega(s) = \dfrac{\dif \mesS}{\dif \mesT}$.
\\
Let $N \in \N^*, \, s_1, \cdots , s_N$ random points of $\Delta \ensS$ drawn with respect to $\mesT$ then an estimator of $\Phi(V(x))$ is given by:
\begin{equation}
\hat{\Phi}_{\mesT} = \dfrac{1}{N}\sum_{i=1}^N \mathbf{1}_{V(x)}(r(s_i))\omega(s_i)
\end{equation}
Furthermore, $\mathbb{E}\hat{\Phi}_{\mesT}^2 = \int_{\Delta \ensS} \mathbf{1}_{V(x)}(r(s))\omega^2(s) \mesT(\dif s)$. If $\mesT \ll \dif s = \dif \thx \dif \thz$, 
\begin{equation}
\mathbb{E}\hat{\Phi}_{\mesT}^2 = \int_{\Delta \ensM} \mathbf{1}_{V(x)}(y)\omega^2(r^{-1}(y))\vert \det{J_r(r^{-1}(y))}\vert^{-1} \dfrac{\dif \mesT}{\dif s}(r^{-1}(y)) \mesM(\dif y)
\end{equation}
\end{propi}
From \ref{th:var} we can provide an approximation of $\mathbb{E}\hat{\Phi}_{\mesT}^2$ and $\sigma_r(\hat{\Phi}_{\mesT}$ in the limit case where $\Phi \ll 1$ and $\mesM(V(x)) \ll H^2$ (and $\mesT \ll \dif s$): 
\begin{equation}
\mathbb{E}\hat{\Phi}_{\mesT}^2 \simeq \omega^2(r^{-1}(x))\vert \det{J_r(r^{-1}(x))}\vert^{-1} \dfrac{\dif \mesT}{\dif s}(r^{-1}(x)) \mesM(V(x))
\end{equation}
\begin{equation}
\sigma_r(\hat{\Phi}_{\mesT}) \simeq \tilde{\sigma}_r(x) \egdef \sqrt{\dfrac{\vert \det{J_r(r^{-1}(x))}\vert}{\dfrac{\dif\mesT}{\dif s}(r^{-1}(x))\mesM(V(x))}}
\end{equation}
We want to design a measure $\mesT$ such as 
\begin{itemize}
\item the sampling according to $\mesT$ is feasible without computation overload (for instance by using the inverse method),
\item the relative variations are sensibly equal through space and the lowest possible
\end{itemize}
The choice of $\mesT$ maximizing $\int \tilde{\sigma_r}^{-1}(x) \dif x$ leads to a $\mesT$ known under closed form: 
\begin{propi}
The maximum of $\int_{\Delta\ensM} \tilde{\sigma}_r^{-1}(x) \dif x$ over measures $\mesT$, $\mesT(\Delta \ensM )=1$, $\mesS \ll \mesT \ll \dif\thx\dif\thz $ is reached for $\dif\mesT = \frac{1}{C_{\mesT}}  \dif \thx \dfrac{1}{H^2}\vert \det{J_r(\thz)}\vert \dif \thz = \frac{1}{C_{\mesT}}  \dif \thx F(\dif \thz)$ with 
\begin{equation}
F(\thz) = \dfrac{1}{\sin^2{\thz_0-\delta\thz}} - \dfrac{1}{\sin^2{\thz}}
\end{equation}
\begin{equation}
C_{\mesT}=\pi\left( \dfrac{1}{\sin^2{\thz_0-\delta\thz}} - \dfrac{1}{\sin^2{\thz_0+\delta\thz}}\right) = \dfrac{\mesM(\Delta \ensM)}{H^2}
\end{equation}
\end{propi}
Thus basically, the relative variation under $\mesM$ is independent of the position and is worth 
\begin{equation}
\sigma_r \simeq \sqrt{ \dfrac{\mesM(V(x))}{N\mesM(\Delta \ensM)}}
\end{equation}
This implies that for $=10^7$ rays sampled from $\mesT$, $\sigma_r \simeq 5\%$ at every point of $\Delta \ensM$.
\section{Results \label{sec:resultats}}
\begin{figure}
\begin{center}
\begin{subfigure}[b]{0.9\textwidth}
                \centering
                \includegraphics[width =\textwidth]{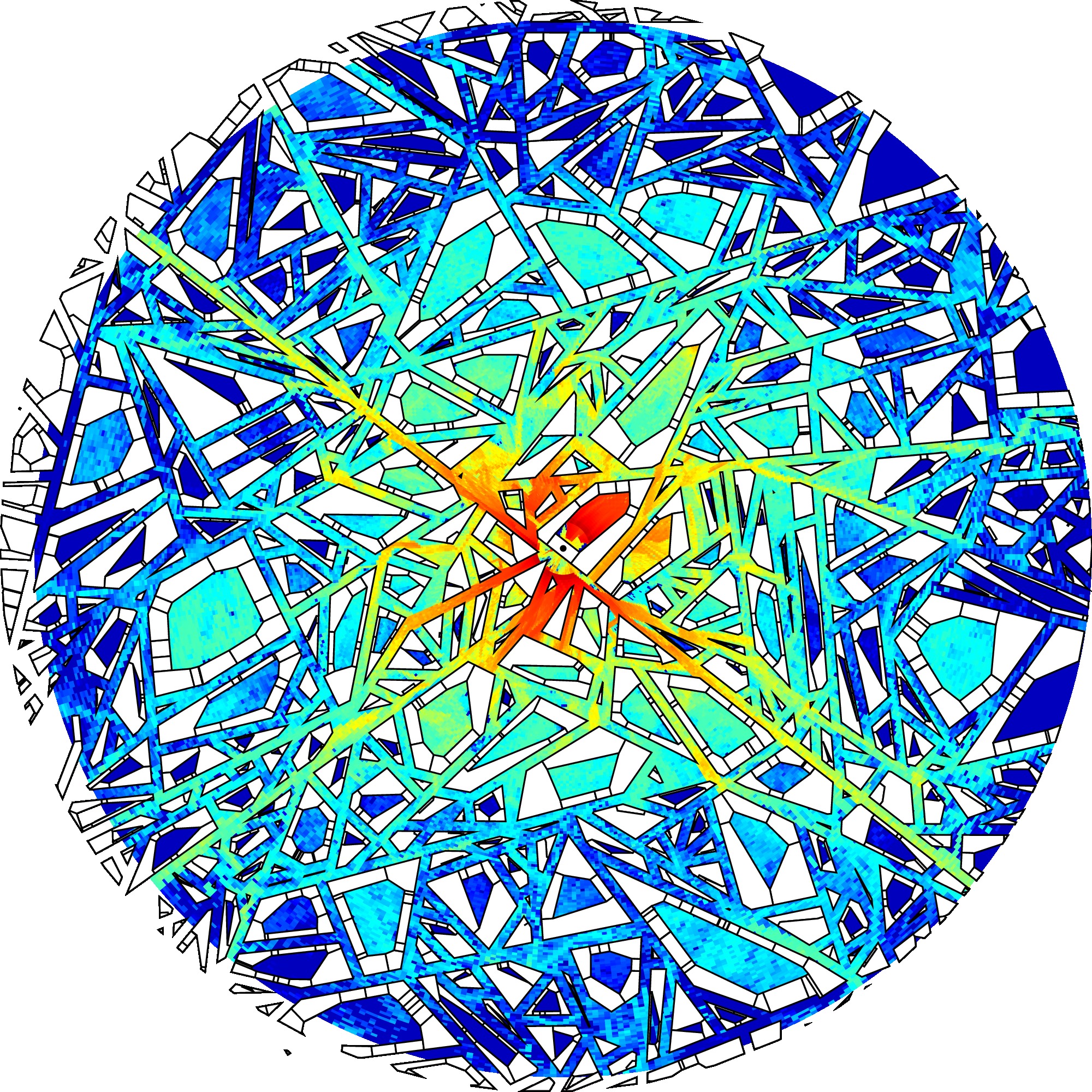}
              
      \end{subfigure}
      
      \end{center}
      \caption{\label{fig:couverture}Estimation of the attenuation map of a random Crack Stit city with typical realistic tunings of \ref{ann:parametres}. (Colors represent the logarithm of the power density). } 
      \end{figure}

The mapping of each point in an horizontal plane the power at this point is called an attenuation map.
\\
In practice, the attenuation map is estimated from pixels surrounding a discrete set of points called evaluation points.
\\  
For instance, Fig.\ref{fig:couverture}, plots the estimated attenuation map for a simulated Crack-STIT city in a disc of radius $1.5$km. The antenna is placed above the roof whose barycenter is the closest to the origin. The evaluation points $x_{j,k}$ have been chosen to form crowns of regularly spaced radius and sectors of constant aperture: $x_{j,k}$ is at the distance $j\delta d$ of the center and forms an angle $k\delta \alpha$ with the $x$-axis. The pixels are arc of crowns surrounding the evaluation point to form a partition of the space.

We notice on this map that: 
\begin{itemize}
\item the power is globally decreasing with the distance
\item some area are shadowed by the buildings 
\item long streets passing next to the antenna behave as wave guides.
\end{itemize}

  To assess the expectancy of power in the street zone, we generate $N$ maps with the same morphological parameters. For each map $i$ we compute the attenuation map $(P_i(d_j, \alpha_k) )_{j,k}$. A pixel that is not at least partially in the street has a power $0$.    
  \begin{equation}
\hat{P}(d_j) = \sum_k\sum_i  P_i(d_j, \alpha_k)\eta_i \dfrac{\delta \alpha}{N \cdot 2\pi}
\end{equation}
where $\eta_i$ the the street area fraction in the simulation window of map $i$ and $\delta \alpha$ is the angular opening of pixels.

Notice that by sum permutation in the above formula, for a particular map, all the pixels in the crown at the distance $d_j$ are averaged, since there are typically $200$ pixels in this crown, the relative error in the estimation of $P_i(d_j)$ is $14$ times smaller than the relative error in the estimation of the power in a particular pixel.

\begin{figure}
\begin{center}
\begin{subfigure}[b]{0.3\textwidth}
                \centering
                \includegraphics[width =\textwidth]{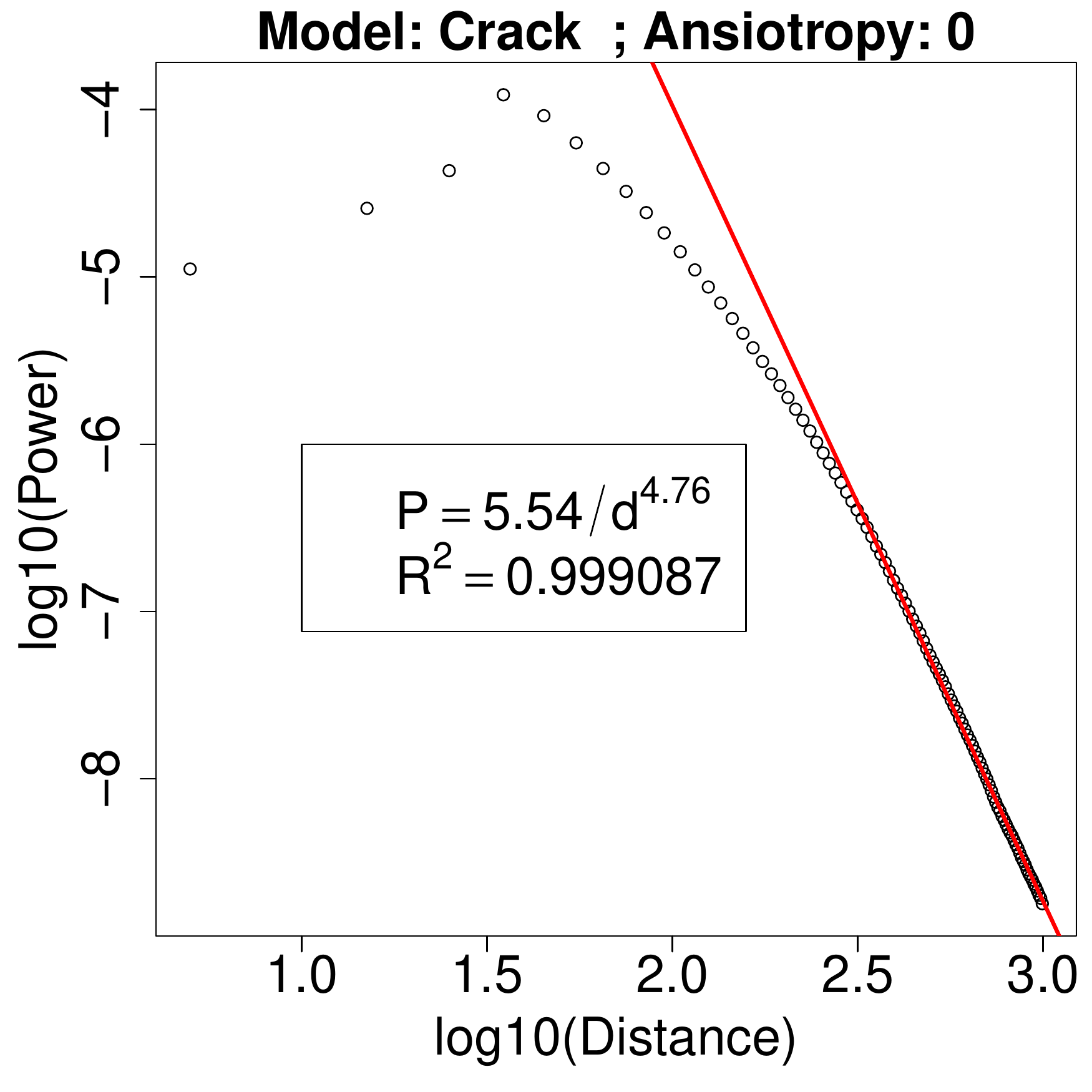}
              
      \end{subfigure}
      \begin{subfigure}[b]{0.3\textwidth}
                \centering
                \includegraphics[width =\textwidth]{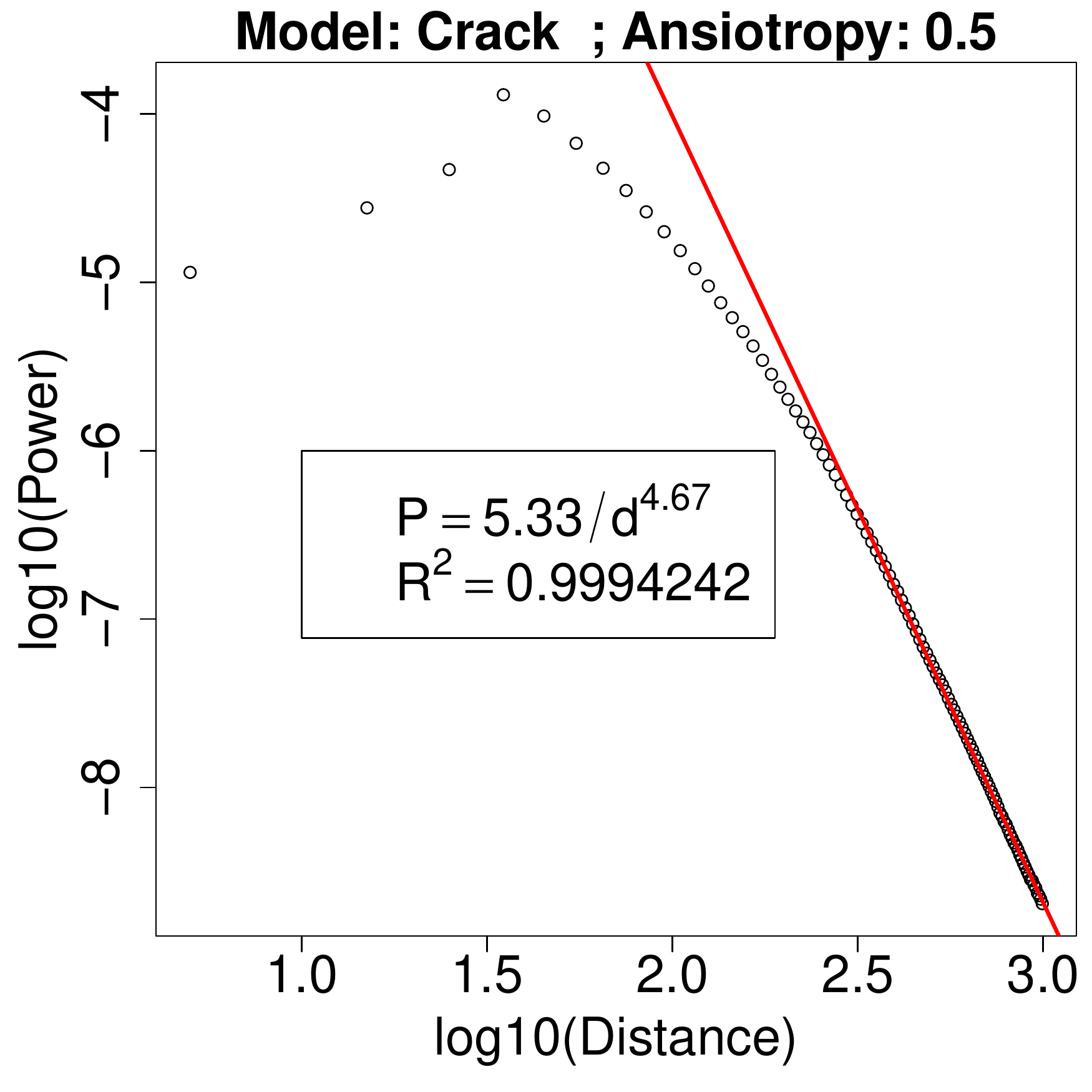}
               
      \end{subfigure}
      \begin{subfigure}[b]{0.3\textwidth}
                \centering
                \includegraphics[width =\textwidth]{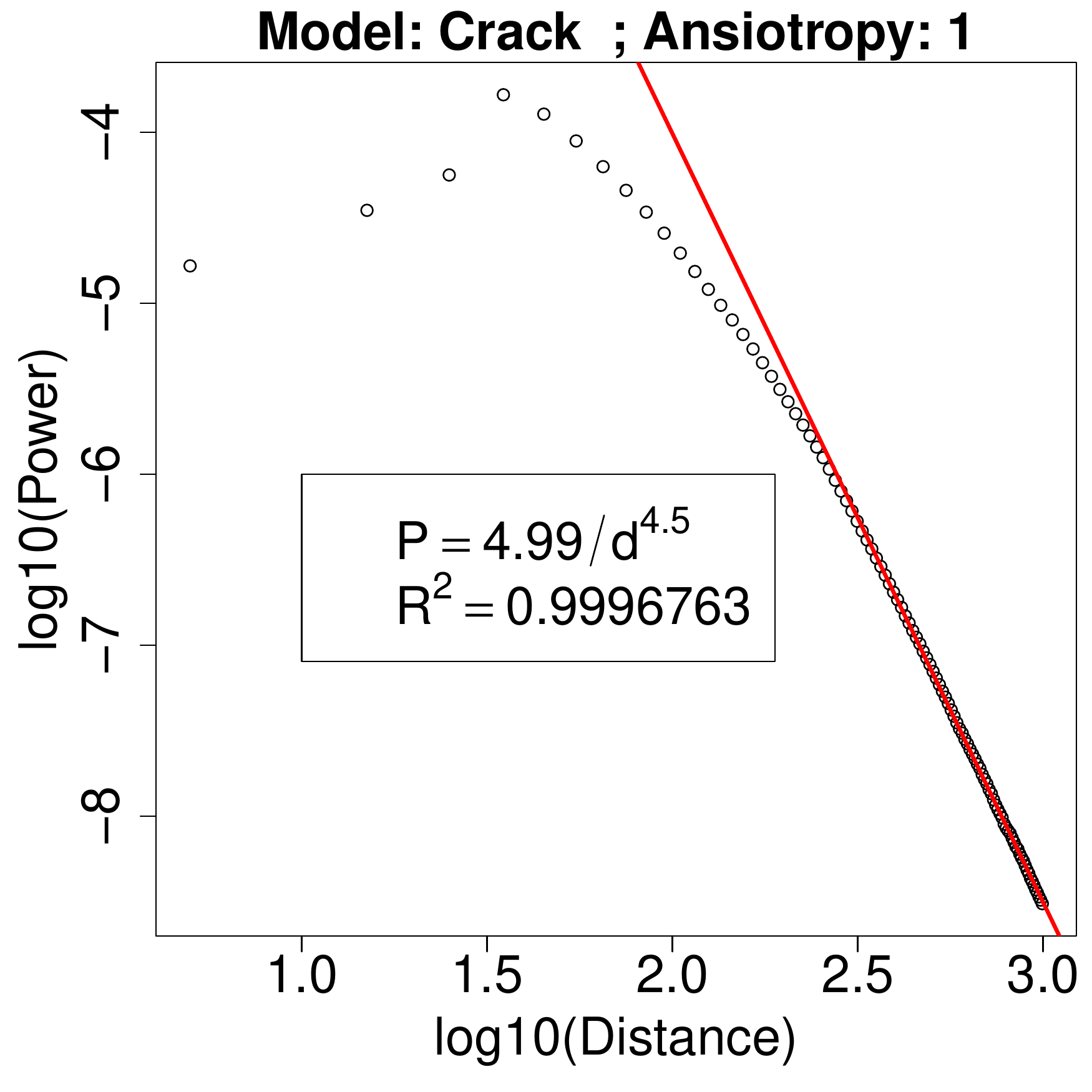}
                
      \end{subfigure}
      \begin{subfigure}[b]{0.3\textwidth}
                \centering
                \includegraphics[width =\textwidth]{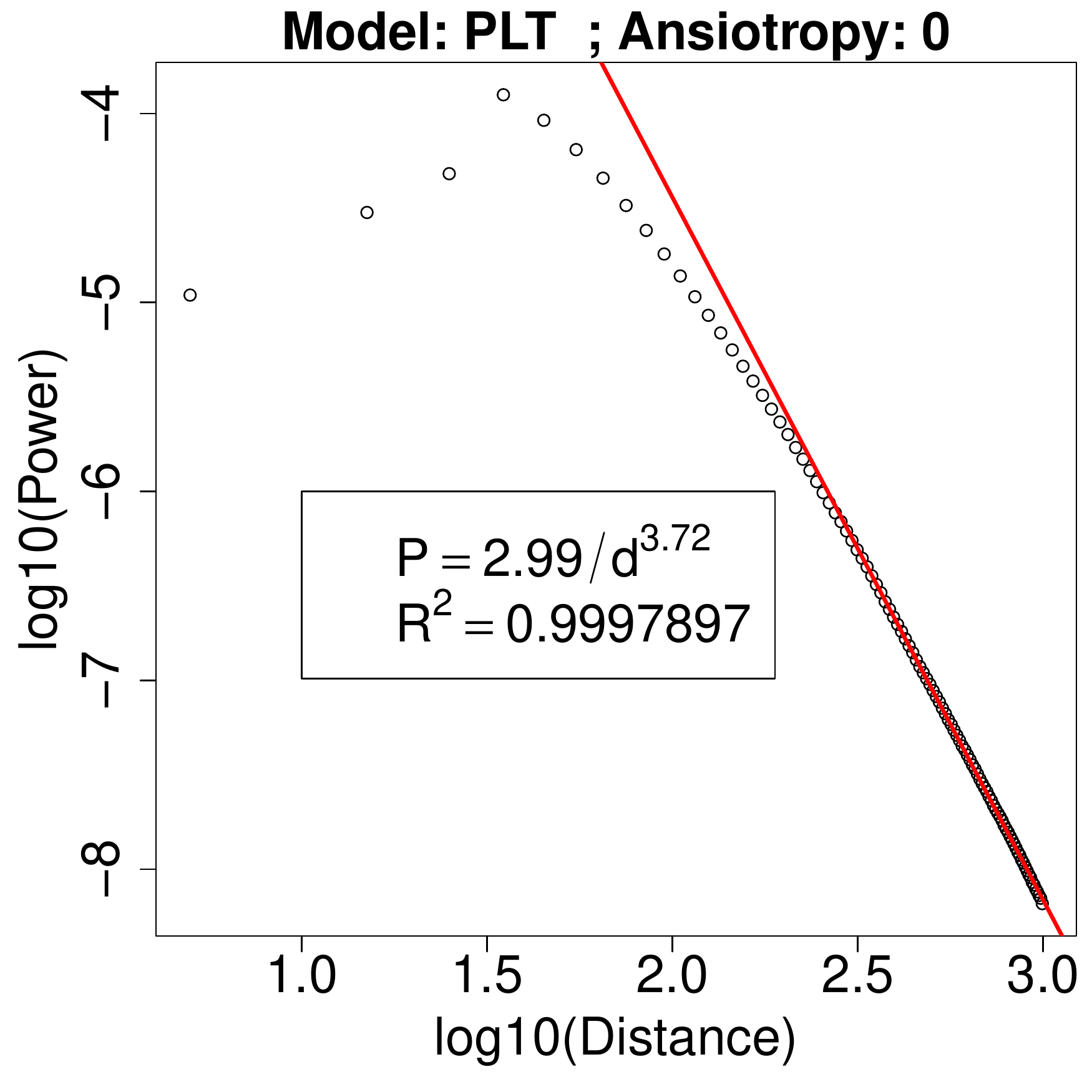}
              
      \end{subfigure}
      \begin{subfigure}[b]{0.3\textwidth}
                \centering
                \includegraphics[width =\textwidth]{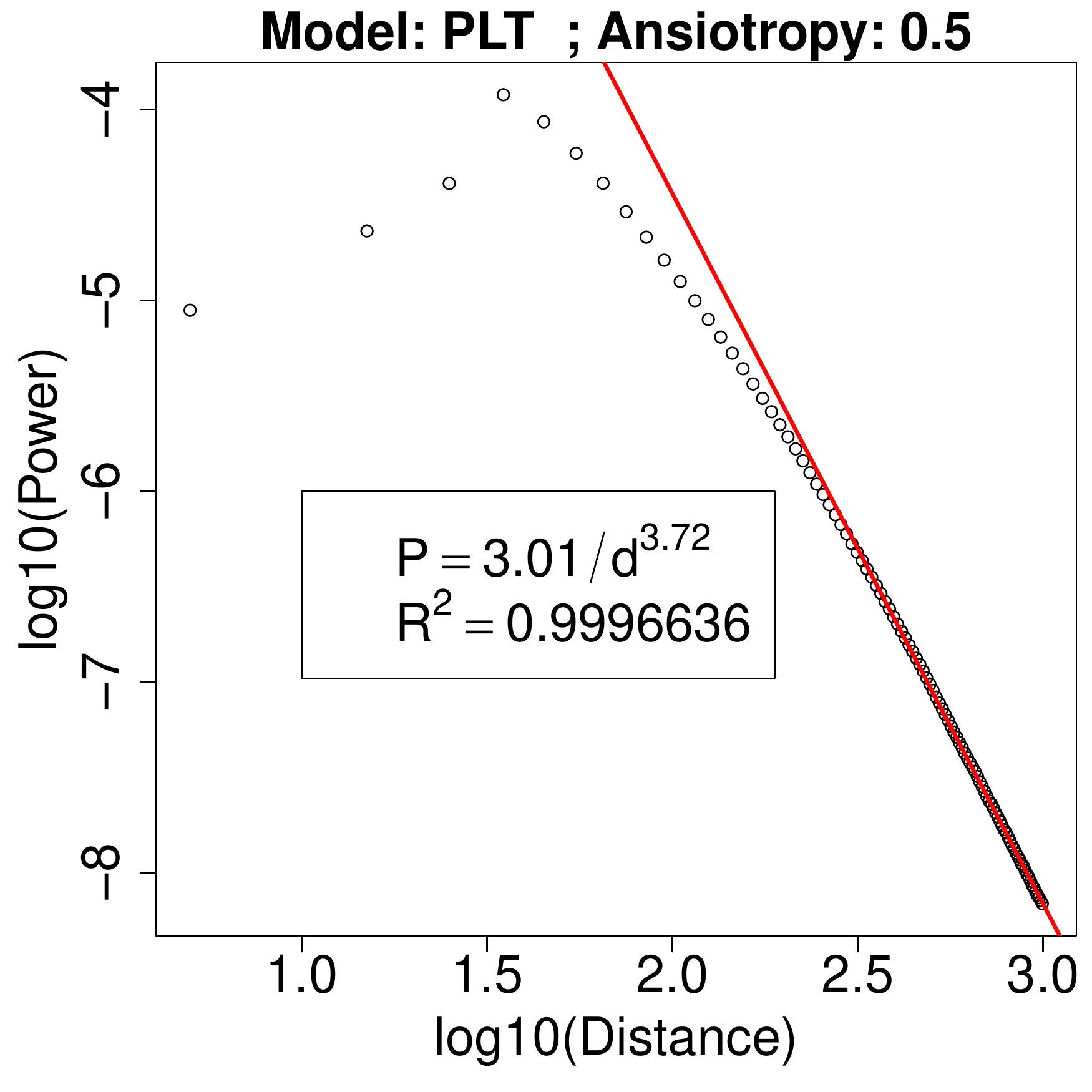}
                
      \end{subfigure}
      \begin{subfigure}[b]{0.3\textwidth}
                \centering
                \includegraphics[width =\textwidth]{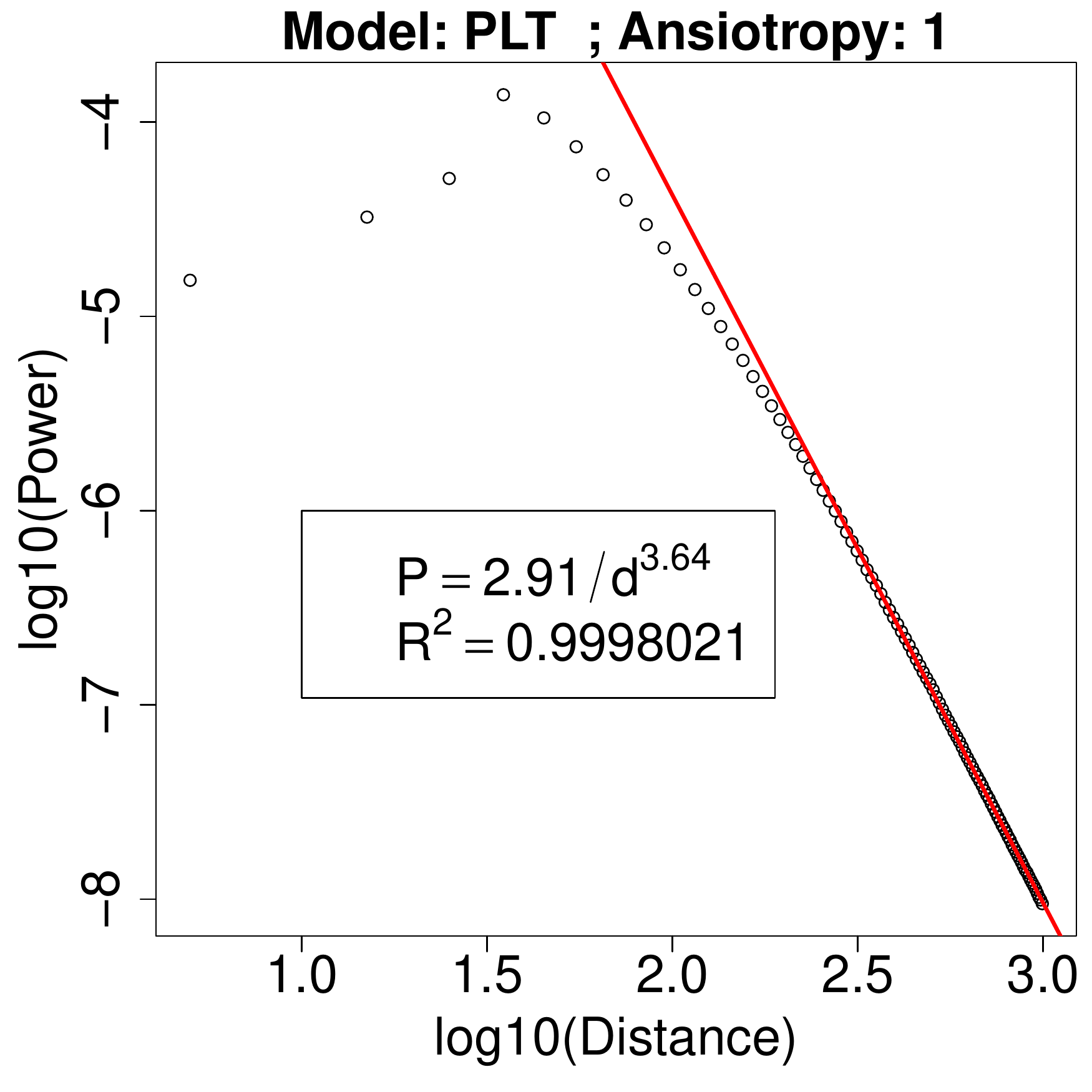}
               
      \end{subfigure}
      \caption{\label{fig:resultats} log-log plots of the mean power received at a distance $d$ for Crack (first row) and PLT (second row). The anisotropy coefficient is 0, 0.5 and 1 from left to right. The fitting by a power function is satisfactory. The attenuation exponent is different for the two classes of tessellations: around $4.6$ for Crack and around $3.7$ for PLT, the latest being a better environment for wave to propagate. The anisotropy does not have a strong influence on the attenuation of the wave.  }
\end{center}
\end{figure}

We have generated for six cases (PLT and Crack STIT both with anisotropy coefficient of $0\, 0.5$ and $1$) $1000$ cities in a disc of radius $R + \Delta R = 1.5$km and estimated the power density in the disc of radius $R$ to avoid side-effects. \\
We obtain estimations of the expected power observed at a certain distance of the source for the different cities' morphology. 
\\
The estimated functions $P(d)$ are fitted by the model $P(d) = A/d^{\alpha}$ (fig.\ref{fig:resultats}). The smallest fits' $R^2-$index is $0.99909$. The order of weight of $\alpha$ is significantly smaller for PLT (around $3.7$) than for Crack (around $4.6$) whereas the anisotropy coefficient has only a slight influence.  At least at a distance equal to twice the typical radius of a serving zone in a cellular network, the path-loss is thus a power function as pointed out empirically in \cite{Sarkar2010}.    
 
\section{Conclusion} 
In this paper we have presented a unified stochastic geometry and ray-tracing framework. This framework has been efficiently implemented in C++ under the name \textit{GeoStat}.
\\
The models of city we have proposed, can be calibrated to fit real world map from simple mean formulas. It mimic main features of real cities like 3D geometry whose propagation properties differs essentially from 2D's ones, fa\c{c}ade alignments along streets that produce wave guide phenomena or organization of buildings into blocks producing wave concentration.   
\\
We have presented a mathematical framework to define and compute properly classical physical quantities in a probabilistic fashion. This framework is general and can be easily adapte to take into account waves interference or propagation schemes different from those of optical geometry.
\\
As an instance of application, we have validated in this paper that the path-loss function is a power function whose attenuation exponent is related to the morphology of the environment at least in the range of distances relevant for telecommunication engineering. A more precise discussion has been presented in \cite{Yu2014}. Other applications such as Sub-Channel distribution has also already been investigated in \cite{Yu2014A} using \textit{GeoStat} framework.
 \appendix

\section{Parameter summary \label{ann:parametres}}
 Here is the summary of the model's parameters and of their notation. To each one is associated a typical numerical value which is used at the different stages of the simulation. We distinguish physical parameters (-) from technical ones (+).  
 
\begin{itemize} \itemsep 0cm \parskip 0cm \parsep 0cm \topsep 0cm
\item Simulation window   
\begin{itemize} \itemsep 0cm \parskip 0cm \parsep 0cm \topsep 0cm
\item[+] $R = 1 \, \text{km} $: radius of the circular window $W$.  
\item[+] $\Delta R = 500 \, \text{m} $: offset radius to avoid side effects. 
\end{itemize}
\item Tessellation 
\begin{itemize} \itemsep0pt \parskip0pt \parsep0pt \topsep0pt
\item $U_2 = 400 \, \text{m}$: mean perimeter of the typical block of houses.
\end{itemize}
\item Buildings
\begin{itemize} \itemsep0pt \parskip0pt \parsep0pt \topsep0pt
\item $l = 10 \, \text{m}$: street thickness.
\item $h = 15 \, \text{m}$: mean height of a building.
\item $b = 10 \, \text{m}$: mean length of a building facade.
\end{itemize}
\item Antenna 
\begin{itemize} \itemsep0pt \parskip0pt \parsep0pt \topsep0pt
\item $P_0 = 40 \, \text{W}$: total power of the source.
\item $f = 2\cdot 10^9 \, \text{Hz}$, $\lambda = 15 \, \text{cm} $: frequency and wave length (GSM).
\item $\theta_0 = \pi / 2 \, \text{rad}$, $\Delta \theta = \pi \, \text{rad}$: horizontal aperture of the source.
\item $\varphi_0 = \pi / 12 \, \text{rad} $, $\Delta \varphi = \pi / 12 \, \text{rad}$: vertical aperture of the source
\end{itemize}
\item Rays
\begin{itemize} \itemsep0pt \parskip0pt \parsep0pt
\item[+] $N = 10^7$: number of rays used in a simulation.
\item $n = \infty$: maximum number of reflections for a ray. 
\item $1/\gamma = 0.5 = -3 \, \text{dB}$: power gain after a reflection. 
\end{itemize}
\item Statistics
\begin{itemize} \itemsep0pt \parskip0pt \parsep0pt
\item[+] $R = 1 \, \text{km} $: radius of the measurement grid.
\item[+] $\delta d = 10 \, \text{m}  $: radial length of a pixel.  
\item[+] $\delta \alpha = 2 \, \text{degrees}$: angular opening of a pixel.
\end{itemize}
\end{itemize} 

\section{Construction of random tessellations \label{ann:tess}} 
We present in this annex the gist of the implementation of the algorithms simulating random tessellations.
\subsection{Data structures and division functions}
\subsubsection{Polygons - Simple Data structures} 
We implement a line, half-line and segment structures that drift from an abstract \textbf{edge} class.
\\
From an edge $e$, we can create the opposite edge $-e$ with the same support and the opposite orientation. 
\\
A \textbf{circular list} $cl$ is a collection of homogeneous objects $cl = (O_1,\cdots, O_n)$ with an iterator $it$ that can be initialized to any $O_i, \, 1 \leq i \leq n$, if at a moment, $it$'s state is $O_j$ then $it.next() = O_{j+1}$ if $j < n$ or $it.next() = O_1$ if $j = n$.
\\ 
A \textbf{convex polygon} is represented by its border, that is to say a circular list of edges $(e_1,\cdots, e_n )$. The integrity of a polygon is ensured by considering edges $e_i$ with the same orientation and sorting them clockwise.  

\subsubsection{Division of a polygon by a line}

The algorithm of division of a polygon $C$ by a line $L$ simply consists in finding intersections  between $L$ with the edges of $C$. 
If there are no intersection, $\mathbf{divise}$ returns $(C, \emptyset)$ or $(\emptyset, C)$ according to the position of $C$ relatively to $L$. \textbf{(Intersection finding step).}

\begin{figure}
\begin{center}

      \begin{subfigure}[b]{0.49\textwidth}
                \centering
                \includegraphics[width =\textwidth]{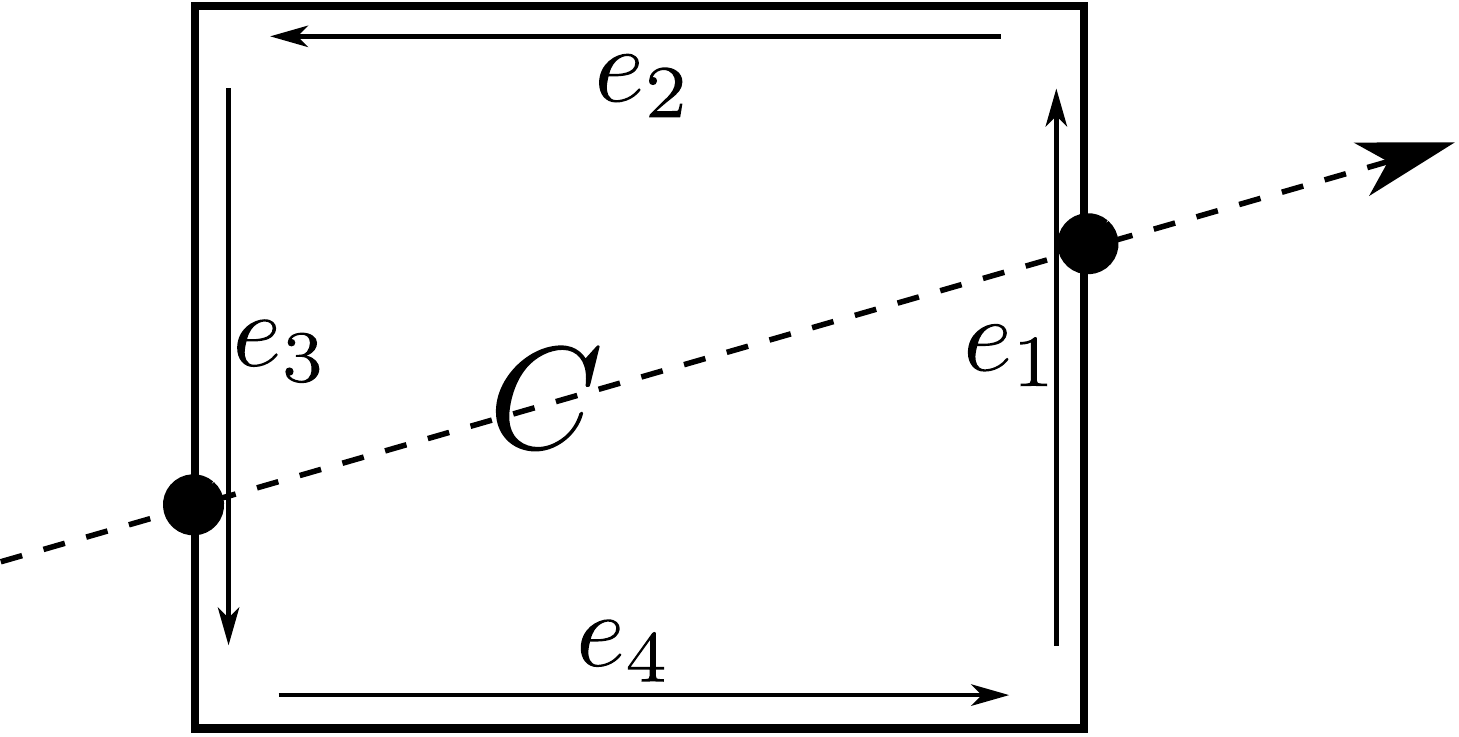}
                
      \end{subfigure}
      \begin{subfigure}[b]{0.49\textwidth}
                \centering
               \includegraphics[width =\textwidth]{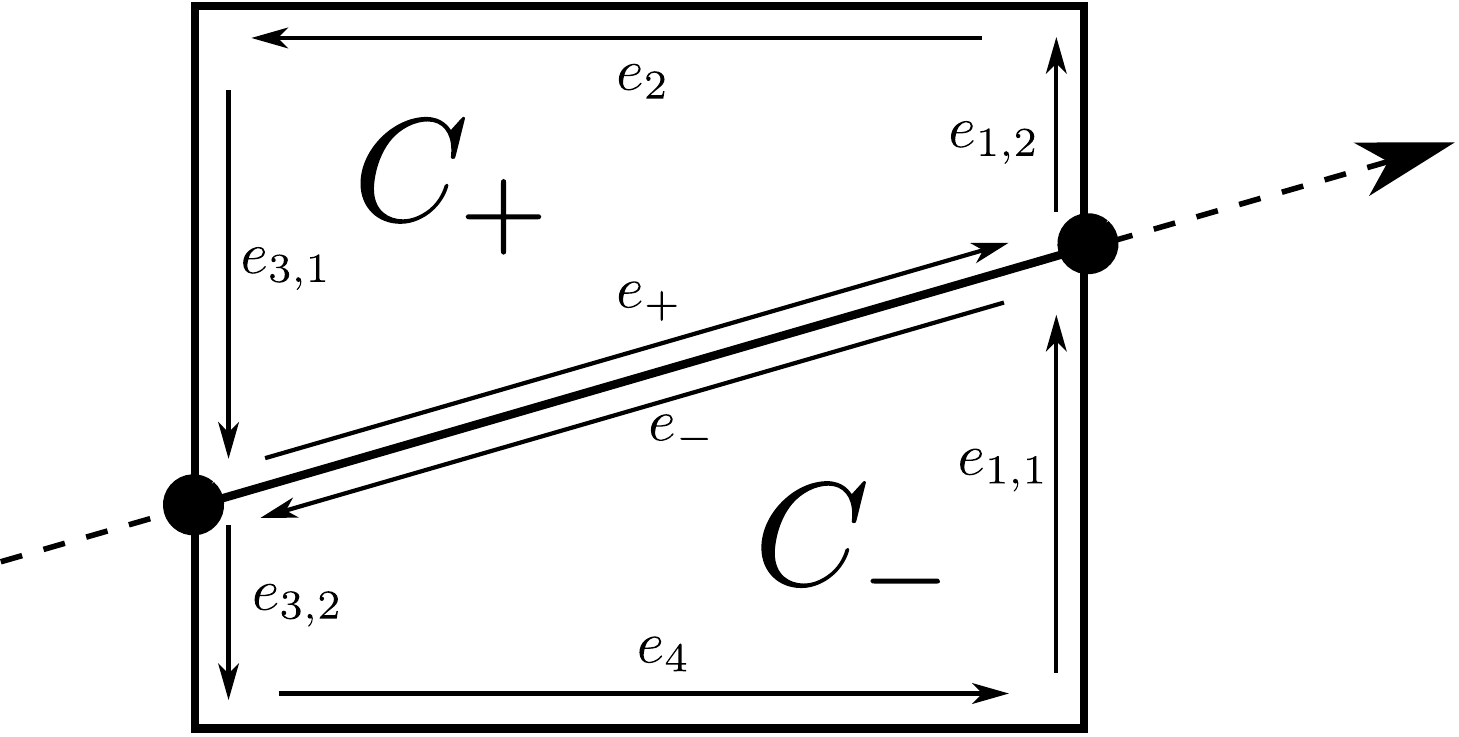}
                
      \end{subfigure}
\end{center}
\caption{\label{fig:division}Division of polygon $C$ by a line into two polygons $C_+$ and $C_-$. In a first step, we look for intersection points $p_1$ and $p_2$. Borders of $C$ are rewritten by introducing theses points and then $C_+$ and $C_-$ can be described by the sequence of their edges.  }

\end{figure}

The most common case is when $L$ intersects two edges (Fig.\ref{fig:division}). 
In this case we write $(f_1, p_1)$ and $(f_2, p_2)$ the intersected edges associated to their intersection points with $L$. Indexes are chosen in such a way that $\overrightarrow{p_1p_2}$ has the same orientation as $L$.  Edges $f_i, \, i=1,2$ give birth to two new edges $f_{i1}$ and $f_{i2}$ with the same orientation as $f_i$. Polygon $C$ that wrote $(e_1, ... , e_k = f_1, ... , e_l = f_2,... e_n)$ (it is always possible to brought back to this representation of the polygon by circular shift of the edges) can be rewritten $(e_1, ..., f_{1,1}, f_{1,2}, .... , f_{2,1}, f_{2,2}, ... e_n)$ . \textbf{(Rewriting step).}

One creates two new segments: $e_+ = \left[p_1,p_2\right]  $ and $e_- = - e_+$. Polygons $C_+$ and $C_-$ resulting from the division can be formally written as : $C_+ =( e_+, f_{1,2}, ... f_{2,1})$ and $C_- =( e_-, f_{2,2},.... f_{1,1})$. \textbf{(Bridge step).}  

The case when a single edge is intersected corresponds to the division of an infinite polygon into two infinite polygons and can be easily written is the same spirit as the previous case.
 
\subsubsection{Tessellation - Complex Data structures}
If the tessellation is seen as a collection of polygons, two problems raise when one wants to divide this tessellation by a line. 
\begin{itemize}
\item It would be necessary to treat the division of each polygon although solely a fraction of the polygons will be actually divided (if the order of magnitude of the number of polygons is $n$, the order of magnitude of divided polygons is $\sqrt{n}$). 
\item The independent divisions of two polygons that share an edge that is cut into a point $p$ produces two new points $p_1$ and $p_2$ theoretically equal but that would not have the same reference in the program execution and that for numerical problems could have a positive distance.
\end{itemize}
To solve these problems, we make a slight change in data structure: edges are replaced by \textbf{edgeContainers} (a structure containing an edge, a pointer to the cell it belongs to and a pointer to the edgeContainer containing its opposed edge) and polygons by \textbf{cells} (a circular list of edgeContainers), a tessellation is then a list of cells:

\begin{itemize} \itemsep0pt \parskip0pt \parsep0pt
\item \textbf{edgeContainer}
\begin{itemize} \itemsep0pt \parskip0pt \parsep0pt
\item Edge * edge
\item edgeContainer * opposite
\item Cell * left
\end{itemize}
\item \textbf{Cell}
\begin{itemize} \itemsep0pt \parskip0pt \parsep0pt
\item Circular list of edgeContainer
\end{itemize}
\item \textbf{Tessellation}
\begin{itemize} \itemsep0pt \parskip0pt \parsep0pt
\item List of cells
\end{itemize}
\end{itemize}

When a cell is divided (Algo.\ref{alg:divtess1}),it is possible to list adjacent cells that are also divided, to update them and to divide them knowing the structures that have been added in the previous cell division (Algo.\ref{alg:divtess2}).

\begin{algorithm}[h!]
\begin{algorithmic}[1]
\caption{\label{alg:divtess1} $\mathbf{division}(C, L, \Xi)$. $C$: cell, $L$: line, $\Xi$: tessellation}
\vspace{0.2cm}
\Statex \textbf{Intersection step}
\vspace{0.2cm}
\State Find edgeContainer $ec_1, ec_2 \in C$ $|$ $L \cap ec_1.\mathbf{edge} = p_1 $: $L \cap ec_2.\mathbf{edge} = p_2 $
\vspace{0.2cm}
\Statex \textbf{Update}
\vspace{0.2cm}
\State $\forall i = 1,2$ cut $ec_i$ by $p_i$ into $ec_{i1}$ and $ec_{i2}$ 
\State Rewrite $C$ in consequence.
\State $\forall i,j$ $ec_{ij}.\mathbf{opposite} \leftarrow  ec_i.\mathbf{opposite}$
\State Create $ \forall i,j \,  \text{minus}ec_{ij} = -  ec_{ij} $
\State Replace $ec_i.\mathbf{opposite}$ by $\text{minus}ec_{i2}, \text{minus}ec_{i1}$ in $ec_i.\mathbf{opposite}.\mathbf{left}$
\State $\forall i,j$ $\text{minus}ec_{ij}.\mathbf{opposite} = ec_{ij}$ 
\vspace{0.2cm}
\Statex \textbf{Bridge step}
\vspace{0.2cm}
\State Exchange $p_1$ and $p_2 \, | \, \langle \overrightarrow{p_1p_2} \vert  L \rangle >0 $ if necessary
\State Create edgeContainer $ec_+$ and $ec_-$ ; $ec_+.\mathbf{opposite} = ec_- $ and conversely. 
\State Create $C_+$ and $C_-$ 
\State $\forall ec \in C_+, \, ec.left = C_+$ idem in $C_-$ 
\vspace{0.2cm}
\Statex \textbf{Tessellation updating}
\vspace{0.2cm}
\State $\Xi = \left( \Xi \backslash C \right) \cup \{ C_+, C_-\} $
\vspace{0.2cm}
\Statex \textbf{Recursive call of the function $\mathbf{division}$}
\vspace{0.2cm}
\State $\mathbf{division}( ec_1.\mathbf{opposite}.\mathbf{left}, L, \Xi , p1 , \text{minus}ec_{1,2})$ (Algo.\ref{alg:divtess2})
\State $\mathbf{division}( ec_2.\mathbf{opposite}.\mathbf{left}, L, \Xi , p2 , \text{minus}ec_{2,2})$ (Algo.\ref{alg:divtess2})
\end{algorithmic}
\end{algorithm}

\begin{algorithm}[h!]
\begin{algorithmic}[1]
\caption{ \label{alg:divtess2} $\mathbf{division}(C, L, \Xi p_0, ec_0 )$  $C$: cell, $L$: line, $\Xi$: tessellation, $p_0$: point in $ec_0$ an edgeContainer.}
\vspace{0.2cm}
\Statex \textbf{Intersection points}
\vspace{0.2cm}
\State $p_1 = p_0$, $ec_1 = ec_0$ 

\State Find $ ec'= ec_2$ $|$ $L \cap ec_2.\mathbf{edge} = p'= p_2 , ec_2 \neq ec_1, \neq ec_1.\mathbf{next}$

\State Cut $ec_2$ by $p_2$ into $ec_{21}$ and $ec_{22}$ 
\vspace{0.2cm}
\Statex \textbf{Updating step}
\vspace{0.2cm}
\State $\forall i$ cut $ec_i$ by $p_i$ into $ec_{i1}$ and $ec_{i2}$ 
\State Rewrite $C$ in consequence.

\State $\forall i,j$ $ec_{ij}.\mathbf{opposite} \leftarrow  ec_i.\mathbf{opposite}$
\State Create $ \forall i,j \,  \text{minus}ec_{ij} = -  ec_{ij} $
\State Remplace $ec_i.\mathbf{opposite}$ by $\text{minus}ec_{i2}, \text{minus}ec_{i1}$ in $ec_i.\mathbf{opposite}.\mathbf{left}$
\State $\forall i,j$ $\text{minus}ec_{ij}.\mathbf{opposite} = ec_{ij}$ 
\vspace{0.2cm}
\Statex \textbf{Bridge step}
\vspace{0.2cm}
\State Exchange $p_1$ and $p_2 \, | \, \langle p_1p_2 \vert L \rangle >0 $ if necessary
\State Create edgeContainer $ec_+$ and $ec_-$ ; $ec_+.\mathbf{opposite} = ec_- $ and conversely. 
\State Create $C_+$ and $C_-$ 
\State $\forall ec \in C_+, \, ec.\textbf{left} = C_+$ idem in $C_-$ 
\vspace{0.2cm}
\Statex \textbf{Tessellation updating}
\vspace{0.2cm}
\State $\Xi = \Xi \backslash C \cup \{ C_+, C_-\} $
\vspace{0.2cm}
\Statex \textbf{Recursive call of the function $\mathbf{division}$}
\vspace{0.2cm}
\State $\mathbf{division}( ec'.\mathbf{opposite}.\mathbf{left}, L, \Xi , p' , \text{minus}ec_{2,2})$
\end{algorithmic}
\end{algorithm}

If there is no intersection, the recursion stops. The division of a tessellation $\Xi$ by a line $L$ goes back to find a cell $C_0$ divided by $L$ and to call $\mathbf{division}(C_0, L, \Xi)$. 

\subsection{Random line simulation}
The tessellations under consideration are constructed from iterated divisions of a polygon by random lines. The following properties allow simulating a random line hitting a given polygon with a probability 1:
\begin{propi}
\label{th:rl1}
If $C$ is a circle with radius $r$ then $\sharp \mathcal{L}_C(.)$ follows a Poisson Law of parameter $\lambda.2.r$ whatever $\mathcal{R}$ is and conditionally to it cuts $C$, a line of the process has a distance to the center of $C$ uniformly distributed. 
\end{propi}
\begin{propi}
\label{th:rl2}
If $W \subset C $ then $\mathcal{L}_W(.) = \mathcal{L}_{W \cap C}(.) =\mathcal{L}_C(.\cap W) $.
\end{propi}
Consequently, it is sufficient to draw a Poisson Line in $C$ and keep only lines that cross $W$ (Algo.\ref{algo:randomline}).
\\
The choice of the smallest circle circumscribed to $W$ permits to minimize line rejections and thus to improve running time.

\begin{algorithm}[h!]
\begin{algorithmic}[1]
\caption{$l \sim \mathcal{L}_W$ Random line in $W$ \label{algo:randomline}}
\State INPUT  $W$ 
\State Inscribe $W$ in a circle of center $0'$ and radius $r$..
\While{}
\State  $r_0 \sim \mathcal{U}_{[-r,r]}$ ; $\alpha \sim \mathcal{R}$
\State Consider the line $l = (r_0, \alpha) + 0'$
\If{$d \cap W \neq \emptyset$}
\Return $l$ 
\State BREAK
\EndIf
\EndWhile 
\end{algorithmic}
\end{algorithm}

\subsection{Algorithms for random line generated processes \label{ann:algo1}}
The division functions (Algo.\ref{alg:divtess1}, \ref{alg:divtess2}) and the random line simulation 
(Algo.\ref{algo:randomline}) based on \ref{th:rl1}, \ref{th:rl2} lead to the simple algorithm \ref{alg:plt1} to simulate the intersection of a PLT with a compact connected window $W$.

\begin{algorithm}[h!]
\begin{algorithmic}[1]
\caption{$\mathbf{\Xi} = \mathbf{PLT}(\lambda, W)$ Poisson Line Tessellation process intersected by $W$  \label{alg:plt1}}
\State INPUT $\lambda \in \mathbb{R}_+$, $W$ 
\State OUTPUT: Tessellation $\mathbf{T}$
\State Tessellation $\mathbf{T}_0 = \{W\}$
\State Inscribe $W$ in a circle f center $0'$ and radius $r$.
\State  $N \sim \mathcal{P}(\lambda.2r)$ ; $n=1$ 
\While{$n \leq N$}
\State  $r_0 \sim \mathcal{U}_{[-r,r]}$ ; $\alpha \sim \mathcal{R}$
\State Consider the line $l = (r_0, \alpha) + 0'$
\If{$d \cap W \neq \emptyset$}
\State $ \mathbf{\Xi}_n = \textbf{division}(\mathbf{\Xi}, l) $
\State $n++$
\Else
\State $ \mathbf{\Xi}_n =  \mathbf{T}_{n-1}$
\EndIf
\EndWhile 
\State $ \mathbf{\Xi} = \mathbf{\Xi}_N$
\end{algorithmic}
\end{algorithm}

The Crack's construction can be made recursively with the function \textbf{division} and a generator of the law $\mathcal{L}_{\omega}$ with $\omega$ a compact set.
\\
To this we define the auxiliary function of evolution of a cell $C$ belonging to a tessellation $\mathbf{\Xi}$ from a time $t>0$: 

\begin{algorithm}[h!]
\begin{algorithmic}[1]
\caption{$\mathbf{evolution}(C , t , \mathbf{\Xi}, \tau, \lambda)$}
\State $\delta t \sim \mathcal{E}(\lambda.\nu(C)$ 
\If{$t+ \delta t < \tau$ }
\State $L \sim \mathcal{R}_C$, $(C_+, C_-) = \textbf{division}(C, L)$
\State$\mathbf{\Xi} = (\mathbf{\Xi} - \{C \} ) \cup \{ C_+ , C_- \}$
\State  $\textbf{evolution}(C_+ , t+ \delta t ,\mathbf{\Xi}, \tau, \lambda)$
\State $\textbf{evolution}(C , t + \delta t,\mathbf{\Xi}, \tau, \lambda)$
\EndIf
\end{algorithmic}
\end{algorithm}
The theorem \ref{prop:crack} rewrites then computationally as: 
\begin{defin}
There exist a stationary, locally finite tessellation whose intersection with a convex and compact window $W$ is the result of
$\mathbf{crack}(W, \tau, \lambda) = \mathbf{evolution}(W , 0 | \{W \},  \tau, \lambda)$. It is called the Crack STIT tessellation.
\end{defin}

\subsection{Running ressources \label{ann:temps}}

To optimize the ray-tracing algorithm, we use a structure that divide recursively the plane into $4^M$ squared regions. Each region contains a list of references to the buildings that intersect them. If $M$ is well chosen, this technique permits to reduce the ray tracing complexity: it becomes almost independent of the number of buildings. It calls for a preprocessing step that can be performed in $\Theta(\log M . N)$.
The table \ref{tab:complexite} sums up the complexity of algorithms at the different steps of the simulation and their average running time on a 32 and 64 bit computer. The most greedy step is the propagation simulation by ray-tracing. The simulation time order of weight is $1.8.10^{-3} N \, \text{s}$ with $N$ the number of rays, other steps being negligible from $4.10^5$ rays. Nonetheless in some particular cases (the antenna is very high or all the buildings have the same height) the algorithm can be improved and the time order of weight decreases to  $2.10^{-4} N \, \text{s}$. 

From a memory point of view, the simulation does not require a lot of resources since rays can be destroyed once their trajectory computation is over and measurement pixels have been updated. 

\begin{table}
\begin{center}
\begin{tabular}{l|c|cc}
\hline
Algorithm & Complexity & 32 bits &  64 bits \\
\hline 
\textbf{Tessellation}  & & & \\ 
PLT &  $(\lambda.r)^{3/2}$ & 0.1s  & 0.03 s \\ 
Crack &  $(\lambda.r)^2$ & 0.08s & 0.02s \\ 
\textbf{Buildings} & $\dfrac{(\lambda .r)^2}{b}$ & 1.2s & 0.25s\\ 
\textbf{Preprocessing}  & $N\cdot \log{M}$ & 5.1s & 1.2s \\ 
\textbf{Propagation} & $N$ & 90 min & 22 min \\ 
\textbf{Statistics} &  $N$ & 1.2s  & 0.4s \\ 
\hline

\end{tabular}
\caption{\label{tab:complexite} Complexity of the different steps of the simulation according to the parameters in \ref{ann:parametres} and the computation time on a computer Linux 32 bits 1.82 GHz, 1Go Ram and a computer Linux 64 bits, 2.93 GHz, 64 Go Ram for typical values of these parameters.  The implement has been performed in C++ and compiled by g++ with -02 option.  }
\end{center}
\end{table}

\bibliographystyle{plain}

\bibliography{biblio}	

\end{document}